\documentclass[12pt]{amsart}
\usepackage{amssymb,amsfonts,amsmath,amsopn,amstext,amscd,color,latexsym,mathabx,mathrsfs,skak,stackrel,xy,url,verbatim}

\input xy
\xyoption{all}

\theoremstyle{remark}

\newtheorem{para}{\bf}[subsection]

\newtheorem{rem}[para]{\bf Remark}
\theoremstyle{definition}

\theoremstyle{plain}

\newtheorem{thm}[para]{Theorem}
\newtheorem{lemma}[para]{Lemma}
\newtheorem{cor}[para]{Corollary}
\newtheorem{prop}[para]{Proposition}

\newenvironment{numequation}{\addtocounter{para}{1}
\begin{equation}}{\end{equation}}

\newcommand{\bbB}{{\mathbb B}}

\newcommand{\bbD}{{\mathbb D}}

\newcommand{\bbG}{{\mathbb G}}

\newcommand{\bbN}{{\mathbb N}}

\newcommand{\bbP}{{\mathbb P}}
\newcommand{\bbQ}{{\mathbb Q}}

\newcommand{\bbX}{{\mathbb X}}

\newcommand{\bbZ}{{\mathbb Z}}

\newcommand{\frg}{{\mathfrak g}}

\newcommand{\frl}{{\mathfrak l}}
\renewcommand{\frm}{{\mathfrak m}}

\newcommand{\fro}{{\mathfrak o}}

\newcommand{\frX}{{\mathfrak X}}

\newcommand{\cD}{{\mathcal D}}

\newcommand{\cI}{{\mathcal I}}

\newcommand{\cO}{{\mathcal O}}

\newcommand{\cR}{{\mathcal R}}

\newcommand{\cT}{{\mathcal T}}

\newcommand{\sD}{{\mathscr D}}

\newcommand{\sH}{{\mathscr H}}

\newcommand{\Qp}{{\mathbb Q_p}}
\newcommand{\Cp}{{\mathbb C_p}}
\newcommand{\Zp}{{\mathbb Z_p}}
\newcommand{\Fp}{{\mathbb F_p}}

\newcommand{\Lie}{{\rm{Lie}}}
\newcommand{\Pf}{{\it Proof. }}

\newcommand{\Spec}{{\rm Spec}}
\newcommand{\Spf}{{\rm Spf}}

\newcommand{\coker}{{\rm coker}}

\newcommand{\lra}{\longrightarrow}
\newcommand{\hra}{\hookrightarrow}
\newcommand{\im}{{\rm im}}

\newcommand{\ra}{\rightarrow}

\newcommand{\sub}{\subset}

\newcommand{\GL}{{\rm GL}}

\newcommand{\Sym}{{\rm Sym}}

\begin{document}

\title{Integral models of $\bbP^1$ and analytic distribution algebras for $\GL_2$}
\author{Deepam Patel}
\address{Institut des Hautes \'Etudes Scientifiques, Le Bois-Marie 35, route de Chartres, 91440 Bures-sur-Yvette, France}
\email{deeppatel1981@gmail.com}
\author{Tobias Schmidt}
\address{Mathematisches Institut, Westf\"alische Wilhelms-Universit\"at M\"unster, Einsteinstr.
62, D-48149 M\"unster, Germany}
\email{toschmid@math.uni-muenster.de}
\author{Matthias Strauch}
\address{Indiana University, Department of Mathematics, Rawles Hall, Bloomington, IN 47405, U.S.A.}
\email{mstrauch@indiana.edu}

\thanks{M. S. would like to acknowledge the support of the National Science Foundation (award DMS-1202303).}

\maketitle

\vskip10pt

\begin{center}
{\it Dedicated to Peter Schneider on the occasion of his sixtieth birthday.}
\end{center}

\vskip30pt

\tableofcontents

\section{Introduction}

The purpose of this paper is to begin the study of connections between arithmetic differential operators on semistable integral and formal models of flag varieties on the one hand and locally analytic distribution algebras of $p$-adic reductive groups on the other hand. Here we only consider the case of the group $\GL_2$ over $\Zp$ and the corresponding flag variety is the projective line $\bbP^1_\Zp$.

\vskip8pt

These investigations are motivated by the wish to study locally
analytic representations of $p$-adic groups geometrically. In  \cite{AW} K. Ardakov and S. Wadsley work with 'crystalline' differential operators (of level zero) on the smooth model of the flag variety of a split reductive group. This is close in spirit
to the classical localization theory of Beilinson-Bernstein
\cite{BB81} and Brylinski-Kashiwara \cite{BK80}. In the paper
\cite{PSS} we have made a first step in merging the localization
theory of Schneider-Stuhler for smooth representations
\cite{SchSt97} with that of \cite{BB81}. A key ingredient is the
embedding, first discovered by V. Berkovich, cf.
\cite{BerkovichBook}, of the building in the non-archimedean
analytic space $X^{\rm an}$ attached to the flag variety $X$ (see also
\cite{RemyThuillierWerner10}). The connection between the building
and $X^{\rm an}$ can also be seen in
terms of formal models for the rigid analytic space $X^{\rm rig}$.
Especially transparent is that relation for formal models of
$\bbP^1$, cf. \cite{BoutotCarayol}. To better understand the
significance of these models for representation theory, and its
relation to distribution algebras, is the starting point for our
work presented here.

\vskip8pt

Regarding distribution algebras, it turns out that the analytic
distribution algebras as considered by M. Emerton in
\cite{EmertonA}, are well suited to be compared to arithmetic
differential operators. Not surprisingly, Emerton has introduced
and studied these rings having Berthelot's theory of arithmetic
differential operators in mind, cf. \cite[sec. 5.2]{EmertonA}.
Arithmetic differential operators on integral smooth models and
their completions have been studied by C. Noot-Huyghe in
\cite{Huyghe97}, \cite{Huyghe09}, \cite{Huyghe98}. In particular,
she proves that these smooth formal models are
$\sD^\dagger$-affine. Here we take up her work in \cite{Huyghe97},
in the special (and easy) case of $\bbP^1$ and show that the ring
of global sections of the arithmetic differential operators is
isomorphic to the analytic distribution algebra
$\cD^{an}(\bbG(0)^\circ)$ of the 'wide open' rigid-analytic group
$\bbG(0)^\circ$ whose $\Cp$-valued points are $\bbG(0)^\circ(\Cp)
= 1+{\rm M}_2(\frm_\Cp)$. Let $\frX$ be the completion of $\bbP^1_\Zp$ along its special fiber. Our first result is

\vskip8pt

{\bf Theorem 1.} (Thm. \ref{thm-0}) {\it There is a canonical isomorphism of (topological) $\Qp$-algebras

$$\cD^{\rm an}(\bbG(0)^\circ)_{\theta_0} \simeq H^0(\frX,\sD^\dagger_{\frX,\bbQ}) \;.$$} \vskip-30pt \qed

\vskip16pt

Here, the subscript $\theta_0$ indicates a central reduction.
The proof of this theorem consists of two parts. Firstly, we identify the analytic distribution algebra $\cD^{\rm an}(\bbG(0)^\circ)$ with the inductive limit (over $m$) of completed 'restricted divided power enveloping algebras' $\widehat{U}(\frg_\Zp)^{(m)}_\bbQ$ (of level $m$) of $\frg_\Zp = \frg\frl_2(\Zp)$.
Secondly, we relate the algebras $U(\frg_\Zp)^{(m)}_\bbQ$ to the global sections $H^0(\bbP^1_\Zp,\cD^{(m)})$ of the sheaf of integral differential operators $\cD^{(m)}$ of level $m$. Much of what we do in this part of the proof (sec. \ref{smooth}) is already contained in \cite{Huyghe97}. We have chosen to redo most of the arguments here, in an entirely explicit manner, because the arguments and techniques will be used later in sections \ref{Xn} and \ref{diff_op_Xn}. 

\vskip8pt

After having obtained theorem 1 we have been informed by C.
Noot-Huyghe that she has proved the general case of this theorem,
for an arbitrary split reductive group and the corresponding smooth  formal model of the flag variety, in an unpublished manuscript. \vskip8pt

Furthermore, we give a description of the analytic distribution algebras $\cD^{\rm an}(\bbG(n)^\circ)$ of rigid-analytic wide open congruence subgroups $\bbG(n)^\circ$.
Their $\Cp$-valued points are given by $\bbG(n)^\circ(\Cp) = 1+p^n{\rm M}_2(\frm_\Cp)$. The description of the distribution algebras is close to that in \cite[sec. 5.2]{EmertonA}, but more suited to the material treated in the
second part of this paper, i.e., sections \ref{Xn} and \ref{diff_op_Xn}. 

\vskip8pt

In these sections we consider
certain semistable integral models $\bbX_n$ of $\bbP^1_\Zp$, and we study the sheaves $\cD_{\bbX_n}^{(m)}$ of logarithmic differential operators of level $m$ on these schemes.
Denote by $H^0(\bbX_n,\cD_{\bbX_n}^{(m)})^\wedge$ the $p$-adic completion of
$H^0(\bbX_n,\cD_{\bbX_n}^{(m)})$, and put $H^0(\bbX_n,\cD_{\bbX_n}^{(m)})^\wedge_\bbQ = H^0(\bbX_n,\cD_{\bbX_n}^{(m)})^\wedge \otimes_\bbZ \bbQ$.
Then we show

\vskip8pt

{\bf Theorem 2.} (Thm. \ref{thm-n}) {\it Given $n \ge 0$ let $n' = \lfloor n\frac{p-1}{p+1} \rfloor$ be the greatest integer less or equal to $n\frac{p-1}{p+1}$. Then we have natural inclusions

$$\cD^{\rm an}(\bbG(n)^\circ)_{\theta_0} \; \hra \; \varinjlim_m H^0(\bbX_n,\cD_{\bbX_n}^{(m)})^\wedge_\bbQ \; \hra \;  \cD^{\rm an}(\bbG(n')^\circ)_{\theta_0} \;.$$} \vskip-30pt \qed

\vskip16pt

Let $\frX_n$ be the formal completion of $\bbX_n$ along its
special fiber. Then there is a canonical injection $\varinjlim_m
H^0(\bbX_n,\cD_{\bbX_n}^{(m)})^\wedge_\bbQ \; \hra \;
H^0(\frX_n,\sD^\dagger_{\frX_n,\bbQ})$. We do not treat here the
question whether this inclusion is in fact an isomorphism. This
problem is related to the question whether the schemes $\bbX_n$
(resp. formal schemes $\frX_n$) are $\cD$-affine, a topic we plan
to discuss in a future paper.

\vskip16pt

{\it Acknowledgements.} The reader will have no difficulty in
recognising the influence of Peter Schneider's work on the ideas
contained in this paper. Over the many years we have spent
together in M\"unster, we have greatly benefited from Peter's
generosity in sharing his ideas with us and guiding us into many
different mathematical worlds. We are grateful for this. It is a
pleasure to dedicate this paper to him on the occasion of his
sixtieth birthday.

\vskip8pt

{\bf Notation.} If $L$ is a field equipped with a non-archimedean absolute value we let $\fro_L$ be its valuation ring and $\frm_{\fro_L} \sub \fro_L$ the maximal ideal of its valuation ring. We let $\bbN = \bbZ_{\ge 0}$ be the set of non-negative integers. If $\nu = (\nu_1,\ldots,\nu_d)$ is a tuple of integers, then we put $|\nu| = \nu_1+\ldots+\nu_d$.

\section{Distribution algebras of wide open congruence subgroups}\label{dist-algs}

\subsection{The group schemes $\bbG(n)$}\label{G(n)}  Let $n \ge 0$ always denote a non-negative integer. Put

$$\bbG(0) = \bbG = \GL_{2,\Zp} = \Spec\left(\Zp\left[a,b,c,d,\frac{1}{\Delta}\right]\right) \;,$$

\vskip8pt

where $\Delta = ad-bc$, and the co-multiplication is the one given by the usual formulas. For $n \ge 1$ let $a_n$, $b_n$, $c_n$, and $d_n$ denote
indeterminates. Define an affine group scheme $\bbG(n)$ over $\Zp$ by setting

$$\cO(\bbG(n)) = \Zp\left[a_n,b_n,c_n,d_n,\frac{1}{\Delta_n}\right] \;, \;\; \mbox{where } \; \Delta_n = (1+p^na_n)(1+p^nd_n)-p^{2n}b_nc_n \;,$$

\vskip8pt

and let the co-multiplication

$$\cO(\bbG(n)) \lra \cO(\bbG(n))\otimes_{\Zp} \cO(\bbG(n)) = \Zp\left[a_n,b_n,c_n,d_n,a_n',b_n',c_n',d_n',\frac{1}{\Delta_n},\frac{1}{\Delta_n'}\right]$$

\vskip8pt

be given by the formulas

$$\begin{array}{lcccccccc}
a_n & \mapsto & a_n &+& a_n' &+& p^na_na_n' &+& p^nb_nc_n' \;,\\
b_n & \mapsto & b_n &+& b_n' &+& p^na_nb_n' &+& p^nb_nd_n' \;,\\
c_n & \mapsto & c_n &+& c_n' &+& p^nc_na_n' &+& p^nd_nc_n' \;,\\
d_n & \mapsto & d_n &+& d_n' &+& p^nd_nd_n' &+& p^nc_nb_n' \;.\\
\end{array}$$

\vskip8pt

These group schemes are connected by homomorphisms $\bbG(n) \ra \bbG(n-1)$ given on the level of algebras as follows:

$$a_{n-1} \mapsto pa_n \;, \;\; b_{n-1} \mapsto pb_n \;, \;\; c_{n-1} \mapsto pc_n \;, \;\; d_{n-1} \mapsto pd_n \;, $$

\vskip8pt

if $n>1$. For $n=1$ we put

$$a \mapsto 1+pa_1 \;, \;\; b \mapsto pb_1 \;, \;\; c \mapsto pc_1 \;, \;\; d \mapsto 1+pd_1 \;. $$

\vskip8pt

For a flat $\Zp$-algebra $R$ the homomorphism $\bbG(n) \ra \bbG(0) = \bbG$ induces an isomorphism of $\bbG(n)(R)$ with a subgroup of $\bbG(R)$, namely

$$\bbG(n)(R)  = \left\{\left(\begin{array}{cc} a & b \\ c & d \end{array}\right) \in \bbG(R) \; \Bigg| \; a-1, b, c, d-1 \in p^nR \;\right\} \;.$$

\vskip8pt

Of course, the preceding formulas defining the group schemes are derived formally from this description by setting $a=1+p^na_n$, $b = p^nb_n$, $c = p^nc_n$, and $d = 1+p^nd_n$.

\subsection{The rigid-analytic groups $\bbG(n)^{{\rm rig}}$ and $\bbG(n)^\circ$}\label{analyticG(n)}  Let $\widehat{\bbG}(n)$ be the completion of $\bbG(n)$ along its special fiber $\bbG(n)_{\Fp}$. This is a formal group scheme over $\Spf(\Zp)$. Its generic fiber in the sense of rigid geometry is an affinoid rigid-analytic group over $\Qp$ which we denote by $\bbG(n)^{\rm rig}$. We have for any completely valued field $L/\Qp$ (whose valuation extends the $p$-adic valuation)

$$\bbG(n)^{\rm rig}(L)  = \left\{\left(\begin{array}{cc} a & b \\ c & d \end{array}\right) \in \bbG(\fro_L) \; \Bigg| \; a-1, b, c, d-1 \in p^n\fro_L \;\right\} \;.$$

\vskip8pt

Furthermore, we let $\widehat{\bbG}(n)^\circ$ be the completion of $\bbG(n)$ in the closed point corresponding to the unit element in $\bbG(n)_{\Fp}$. This is a formal group scheme over $\Spf(\Zp)$ (not of topologically finite type). Its generic fiber in the sense of Berthelot, cf. \cite[sec. 7.1]{deJongCrystalline}, is a so-called 'wide open' rigid-analytic group over $\Qp$ which we denote by $\bbG(n)^\circ$. We have for any completely valued field $L/\Qp$ (whose valuation extends the $p$-adic valuation)

$$\bbG(n)^\circ(L)  = \left\{\left(\begin{array}{cc} a & b \\ c & d \end{array}\right) \in \bbG(\fro_L) \; \Bigg| \; a-1, b, c, d-1 \in p^n\frm_{\fro_L} \;\right\} \;.$$

\vskip8pt

The remainder of this section is inspired by M. Emerton's paper \cite{EmertonA}, especially sec. 5.

\subsection{The analytic distribution algebra of $\bbG(0)^\circ$}\label{analytic_dist_alg_G(0)} Our goal in this subsection is to give a description of

$$\cD^{\rm an}(\bbG(0)^\circ) \stackrel[\rm df]{}{=} \cO(\bbG(0)^\circ)'_b$$

\vskip8pt

in terms of 'divided power enveloping algebras' which is analogous to \cite[5.2.6]{EmertonA}. However, the discussion in \cite[sec. 5.2]{EmertonA} does not apply here because the exponential function for the group $\GL_2(\Qp)$ does not map a lattice in

$$\frg \stackrel[\rm df]{}{=} \Lie(\GL_2(\Qp))$$

\vskip8pt

bijectively onto $\GL_2(\Zp)$. Nevertheless, it is possible to also treat $\bbG(0)^\circ$ by making use of the 'Kostant $\bbZ$-form' of the enveloping algebra $U(\frg)$, cf. \cite{Kostant66}. Set

$$e = \left(\begin{array}{cc} 0 & 1 \\ 0 & 0 \end{array}\right) \;, \;\;  h_1 = \left(\begin{array}{cc} 1 & 0 \\ 0 & 0 \end{array}\right) \;, \;\;  h_2 = \left(\begin{array}{cc} 0 & 0 \\ 0 & 1 \end{array}\right) \;, \;\;  f = \left(\begin{array}{cc} 0 & 0 \\ 1 & 0 \end{array}\right) \;,$$

\vskip8pt

and put

$$\frg_\Zp \stackrel[\rm df]{}{=} {\rm M_2}(\Zp) = \Zp e \oplus \Zp h_1 \oplus \Zp h_2 \oplus \Zp f \;.$$

\vskip8pt

For integers $m, n \in \bbN$ define

$$q^{(m)}_n \stackrel[\rm df]{}{=} \left\lfloor \frac{n}{p^m} \right\rfloor \;,$$

\vskip8pt

that is, the greatest integer less or equal to $\frac{n}{p^m}$.
For fixed $m$ we then denote by $U(\frg_\Zp)^{(m)}$ the $\Zp$-submodule of $U(\frg)$ generated by the elements

\begin{numequation}\label{generators}q^{(m)}_{\nu_1}! \frac{e^{\nu_1}}{\nu_1!}  \cdot q^{(m)}_{\nu_2}! {h_1 \choose \nu_2} \cdot q^{(m)}_{\nu_3}! {h_2 \choose \nu_3} \cdot q^{(m)}_{\nu_4}! \frac{f^{\nu_4}}{\nu_4!} \;.
\end{numequation}

\begin{lemma}\label{submod_is_subalg} $U(\frg_\Zp)^{(m)}$ is a $\Zp$-subalgebra of $U(\frg)$.
\end{lemma}

\Pf This is contained in \cite[Prop. 2.3.1]{Montagnon} and the remark before \cite[Lemme 2.3.3]{Montagnon}, namely that $\cD^{(m)}_{X,n}$ has a basis given by the operators $\partial_{<\underline{k}>}$, $|\underline{k}| \le n$. Note also the description of $\partial_{<\underline{k}>}$ given in part (c) of that lemma. \qed

\vskip8pt

We now let $\widehat{U}(\frg_\Zp)^{(m)}$ be the $p$-adic completion of $U(\frg_\Zp)^{(m)}$. Explicitly, its elements can be written as

$$\sum_{\nu = (\nu_1,\nu_2,\nu_3,\nu_4) \in \bbN^4} \gamma_\nu \cdot q^{(m)}_{\nu_1}! \frac{e^{\nu_1}}{\nu_1!}  \cdot q^{(m)}_{\nu_2}! {h_1 \choose \nu_2} \cdot q^{(m)}_{\nu_3}! {h_2 \choose \nu_3} \cdot q^{(m)}_{\nu_4}! \frac{f^{\nu_4}}{\nu_4!} \;,$$

\vskip8pt

where $\gamma_\nu \in \Zp$ and $|\gamma_\nu| \ra 0$ as $|\nu| \ra \infty$. Furthermore, we put

$$\widehat{U}(\frg_\Zp)^{(m)}_\bbQ \stackrel[\rm df]{}{=} \widehat{U}(\frg_\Zp)^{(m)} \otimes_\bbZ \bbQ \;.$$

\vskip8pt

We consider the unique $\Qp$-algebra homomorphism $U(\frg) \ra \cD^{\rm an}(\bbG(0)^\circ)$ which sends $X \in \frg$ to
the linear form

$$f \mapsto \frac{d}{dt} f(e^{tX})\Big|_{t=0} \;.$$

\vskip8pt

Here we follow the same convention as in \cite[sec. 5]{EmertonA} in that we consider the right regular action of a group on its ring of functions.

\begin{prop}\label{dist_alg_G(0)} The map $U(\frg) \ra \cD^{\rm an}(\bbG(0)^\circ)$ just defined extends continuously to $\widehat{U}(\frg_\Zp)^{(m)}$. The family of these maps, for various $m$, induces a canonical isomorphism of topological $\Qp$-algebras

$$\varinjlim_m \widehat{U}(\frg_\Zp)^{(m)}_\bbQ \stackrel{\simeq}{\lra} \cD^{\rm an}(\bbG(0)^\circ) \;.$$

\vskip8pt
\end{prop}

\Pf The affine algebra of the formal group scheme $\widehat{\bbG}(0)^\circ$ is the completion of the ring
$\Zp\left[a,b,c,d,\frac{1}{\Delta}\right]$ with respect to the ideal $I = (p,a-1,b,c,d-1)$. (We write here $a$ instead of $a_0$, $b$ instead of $b_0$, etc.) Hence

$$\cO\left(\widehat{\bbG}(0)^\circ\right) = \Zp[[a-1,b,c,d-1]] \;.$$

\vskip8pt

For the ring of global functions of $\bbG(0)^\circ$ we then have, algebraically and topologically,

$$\cO\left( \bbG(0)^\circ \right) = \varprojlim_{r<1} \cO\left( \bbG(0)_r \right) \;,$$

where

$$\cO\left( \bbG(0)_r \right) = \left\{ \sum \xi_\mu (a-1)^{\mu_1}b^{\mu_2}c^{\mu_3}(d-1)^{\mu_4} \; \Bigg| \; |\xi_\mu|r^{|\mu|} \ra 0 \mbox{ as } |\mu| \ra \infty \right\} \;.$$

\vskip8pt

It is easily checked that

$$\left[\frac{e^{\nu_1}}{\nu_1!} {h_1 \choose \nu_2} {h_2 \choose \nu_3} \frac{f^{\nu_4}}{\nu_4!}\right] . \left[(a-1)^{\mu_1}b^{\mu_2}c^{\mu_3}(d-1)^{\mu_4} \right] = \left\{
\begin{array}{lcl}
1 & , & \nu = \mu \\
0 & , & \nu \neq \mu \end{array} \right. \;.$$

\vskip8pt

We thus find that $\cD^{\rm an}(\bbG(0)^\circ)$ consists of sums

$$\sum_{\nu = (\nu_1,\nu_2,\nu_3,\nu_4) \in \bbN^4} \gamma_\nu \frac{e^{\nu_1}}{\nu_1!} {h_1 \choose \nu_2} {h_2 \choose \nu_3} \frac{f^{\nu_4}}{\nu_4!} \;,$$

\vskip8pt

which have the property that there is $R>1$ for which $|\gamma_\nu|R^{|\nu|} \ra 0$ as $|\nu| \ra \infty$. The rest of the proof is as in \cite[5.2.6]{EmertonA}. Because

$$v_p\left(q^{(m)}_{\nu_1}! q^{(m)}_{\nu_2}! q^{(m)}_{\nu_3}!  q^{(m)}_{\nu_4}!\right) $$

\vskip8pt

is asymptotic to

$$\frac{\nu_1+\nu_2+\nu_3+\nu_4}{(p-1)p^m} \;\; \mbox{ as } \;\; |\nu| \ra \infty \;,$$

\vskip8pt

it follows that $\widehat{U}(\frg_\Zp)^{(m)}_\bbQ$ embeds into $\cD^{\rm an}(\bbG(0)^\circ)$. Furthermore, the inductive limit of the spaces $\cO(\bbG(0)_r)'_b$, for $r \uparrow 1$, is equal to the the inductive limit of the rings $\widehat{U}(\frg_\Zp)^{(m)}_\bbQ$, as $m \ra \infty$. \qed

\vskip8pt

\begin{rem} The Kostant $\bbZ$-form of $U(\frg)$ is nothing else than the distribution algebra ${\rm Dist}(\GL_{2,\Zp})$ of the group scheme $\GL_{2,\Zp}$ as defined in \cite[I.7]{Jantzen}, cf. \cite[II.1.12]{Jantzen} for the explicit relation between the Kostant $\bbZ$-form and the distribution algebra. One can then use the very definition of
the distribution algebra in \cite[I.7]{Jantzen} to give an intrinsic proof of \ref{dist_alg_G(0)} which should generalize to any split reductive group scheme over $\Zp$. \qed
\end{rem}

\subsection{The analytic distribution algebra of $\bbG(n)^\circ$ for $n \ge 1$}\label{dist_alg_G(n)_defn} In this subsection we derive a description of $\cD^{\rm an}(\bbG(n)^\circ) \stackrel[\rm df]{}{=} \cO(\bbG(n)^\circ)'_b$, for $n \ge 1$, from the decription in \ref{dist_alg_G(0)}. The open embedding of rigid spaces $\bbG(n)^\circ \hra \bbG(0)^\circ$ induces a restriction map on spaces of functions $\cO(\bbG(0)^\circ) \ra \cO(\bbG(n)^\circ)$ which has dense image. Taking the continuous dual spaces gives hence an injection

$$\cD^{\rm an}(\bbG(n)^\circ) \hra \cD^{\rm an}(\bbG(0)^\circ) \;.$$

\vskip8pt

We will describe the left hand side as a subalgebra of the right hand side. To this end, let $U(p^n\frg_\Zp)^{(m)}$ be the $\Zp$-submodule of $U(\frg)$ generated by the elements

\begin{numequation}\label{generators_n}
q^{(m)}_{\nu_1}! \frac{(p^ne)^{\nu_1}}{\nu_1!}  \cdot q^{(m)}_{\nu_2}! p^{n\nu_2}{h_1 \choose \nu_2} \cdot q^{(m)}_{\nu_3}! p^{n\nu_3} {h_2 \choose \nu_3} \cdot q^{(m)}_{\nu_4}! \frac{(p^nf)^{\nu_4}}{\nu_4!} \;.
\end{numequation}

As before, we find that $U(p^n\frg_\Zp)^{(m)}$ is a $\Zp$-subalgebra of $U(\frg)$, and we let $\widehat{U}(p^n\frg_\Zp)^{(m)}$ denote its $p$-adic completion.

\vskip8pt

\begin{rem}\label{caution} We caution the reader that $U(p^n\frg_\Zp)^{(m)}$ and $\widehat{U}(p^n\frg_\Zp)^{(m)}$ are merely notations. That is, these rings are not what one would get by formally replacing (the basis of) $\frg_\Zp$ by (the basis of) $p^n\frg_\Zp$ in the definition of
$U(\frg_\Zp)^{(m)}$. The reason is that, obviously,

$${p^n h_i \choose \nu} \; \neq p^{n\nu} {h_i \choose \nu} \;,$$

\vskip8pt

if $\nu > 1$. It is the term on right which one has to work with here, not the term on the left. \qed
\end{rem}

The algebra homomorphism $U(\frg) \ra \cD^{\rm an}(\bbG(0)^\circ)$ defined right before \ref{dist_alg_G(0)} obviously factors as $U(\frg) \ra \cD^{\rm an}(\bbG(n)^\circ) \ra \cD^{\rm an}(\bbG(0)^\circ)$.

\begin{prop}\label{dist_alg_G(n)} The map $U(\frg) \ra \cD^{\rm an}(\bbG(n)^\circ)$ extends continuously to $\widehat{U}(p^n\frg_\Zp)^{(m)}$ and there is a canonical isomorphism of topological $\Qp$-algebras

$$\varinjlim_m \widehat{U}(p^n\frg_\Zp)^{(m)}_\bbQ \stackrel{\simeq}{\lra} \cD^{\rm an}(\bbG(n)^\circ) \;.$$

\vskip8pt
\end{prop}

\Pf We proceed here as in the proof of \ref{dist_alg_G(0)}. The affine algebra of the formal group scheme $\widehat{\bbG}(n)^\circ$ is $\Zp[[a_n,b_n,c_n,d_n]]$ and the coordinates $a_n, b_n, c_n, d_n$ on $\bbG(n)^\circ$ are related to the coordinates $a,b,c,d$ on $\bbG(0)^\circ$ by

$$a_n = \frac{1}{p^n}(a-1)\;, \hskip8pt
b_n = \frac{1}{p^n}b\;, \hskip8pt
c_n = \frac{1}{p^n}c\;, \hskip8pt
d_n = \frac{1}{p^n}(d-1)\;.$$

\vskip8pt

>From the proof of \ref{dist_alg_G(0)} we get

\vskip8pt

$\left[\frac{(p^ne)^{\nu_1}}{\nu_1!} p^{n\nu_2}{h_1 \choose \nu_2} p^{n\nu_3}{h_2 \choose \nu_3} \frac{(p^nf)^{\nu_4}}{\nu_4!}\right] . \left[\left(\frac{a-1}{p^n}\right)^{\mu_1}\left(\frac{b}{p^n}\right)^{\mu_2}\left(\frac{c}{p^n}\right)^{\mu_3}\left(\frac{d-1}{p^n}\right)^{\mu_4} \right] = \left\{
\begin{array}{lcl}
1 & , & \nu = \mu \\
0 & , & \nu \neq \mu \end{array} \right. \;.$

\vskip8pt

And the remainder of the proof is along the same lines as in \ref{dist_alg_G(0)}. \qed

\begin{rem}\label{comp_to_Emerton} For $n \ge 1$ ($n \ge 2$ if $p=2$) the group $\bbG(n)(\Zp) = 1+p^n{\rm M}_2(\Zp)$ is uniform pro-$p$ and its integral Lie algebra in the sense of \cite[sec. 9]{DDMS} is $p^n\frg_\Zp$ when considered as a $\Zp$-submodule of $\frg$. We can thus apply \cite[sec. 5.2]{EmertonA} to get a description of $\cD^{\rm an}(\bbG(n)^\circ)$ in terms of divided power enveloping algebras. The relation between the two descriptions is as follows.

\vskip8pt

In \cite{EmertonA}, $\bbG(n)^\circ$ is identified with the rigid-analytic four-dimensional wide open polydisc $(\bbB^\circ)^4$ via the 'coordinates of the second kind'

$$(t_1,t_2,t_3,t_4) \mapsto \exp(t_1p^ne) \exp(t_2p^nh_1) \exp(t_3p^nh_2) \exp(t_4p^nf) \;.$$

\vskip8pt

Functions $\cO(\bbG(n)^\circ)$ are then considered as functions on $(\bbB^\circ)^4$ via pull-back. Using this identification, we consider elements in $U(\frg)$ as differential operators on $\cO\left((\bbB^\circ)^4\right)$. \cite[5.2.6]{EmertonA} then tells us that $\cD^{\rm an}(\bbG(n)^\circ)$ is the inductive limit of rings

$$\begin{array}{l}
\cD^{\rm an}(\bbG(n)^\circ)^{(m)}\\
\\
\stackrel[\rm df]{}{=} \left\{ \sum_\nu \gamma_\nu \frac{q^{(m)}_{\nu_1}! q^{(m)}_{\nu_2}! q^{(m)}_{\nu_3}! q^{(m)}_{\nu_4}!}{\nu_1! \nu_2! \nu_3! \nu_4!} (p^ne)^{\nu_1} (p^nh_1)^{\nu_2} (p^nh_2)^{\nu_3} (p^nf)^{\nu_4} \; \Bigg| \; |\gamma_\nu| \ra 0 \mbox{ as } |\nu| \ra 0 \right\}
\;. \end{array}$$

\vskip8pt

The relation of these rings to the rings $\widehat{U}(p^n\frg_\Zp)^{(m)}_\Qp$ follows immediately from the elementary

\begin{prop} Suppose $n \ge 1$ ($n \ge 2$ if $p=2$), and let $T$ be an indeterminate. For all $\nu \ge 0$, if one writes the
polynomial $p^{n\nu}{T \choose \nu}$ as

$$\sum_{j=1}^\nu c_{\nu,j} \frac{(p^nT)^j}{j!} \;,$$

\vskip8pt

the coefficients $c_{\nu,j}$ are in $\Zp$.

\end{prop}

\Pf Let $z$ be another indeterminate and consider the formal power series

$$\sum_{\nu \ge 0} p^{n\nu}{T \choose \nu}z^\nu \;.$$

\vskip8pt

This is equal to $(1+p^nz)^T = \exp(T\log(1+p^nz))$. Under the assumption $n \ge 1$ ($n \ge 2$ if $p=2$), one can write $\log(1+p^nz) = p^nzf(z)$ with a power series $f(z) \in \Zp[[z]]$. Hence

$$\exp(T\log(1+p^nz)) = \sum_{j \ge 0} (zf(z))^j \frac{(p^nT)^j}{j!} \;.$$

\vskip8pt

Now compare the coefficients of $z^\nu$ on both sides. \qed

\end{rem}

\section{Arithmetic differential operators on the smooth formal model}\label{smooth}

\subsection{Differential operators with divided powers} We consider $\bbX \stackrel[\rm df]{}{=} \bbP^1_\Zp$ as being glued together from the affine lines

$$U_x = \Spec(\Zp[x]) \;\; \mbox{ and } \;\; U_y = \Spec(\Zp[y])$$

\vskip8pt

along the open subsets $\Spec(\Zp[x,\frac{1}{x}])$ and $\Spec(\Zp[y,\frac{1}{y}])$ according to the relation $xy=1$. The formulas

$$x.\left(\begin{array}{cc} a&b\\c&d \end{array}\right) = \frac{b+dx}{a+cx} \;\;, \hskip10pt y.\left(\begin{array}{cc} a&b\\c&d \end{array}\right) = \frac{ay+c}{by+d} \;,$$

\vskip8pt

describe a {\it right} action of $\bbG = \GL_{2,\Zp} = \Spec\left(\Zp\left[a,b,c,d,\frac{1}{\Delta}\right]\right)$ on $\bbX$. Put $\partial_x = \frac{d}{dx}$ and $\partial_y = \frac{d}{dy}$. These differential operators satisfy the relations

$$\partial_x = -y^2 \partial_y\;, \;\; x\partial_x = -y \partial_y\;, \;\; x^2\partial_x = -\partial_y \;.$$

\vskip8pt

Denote by $\cT_\bbX$ the tangent sheaf of $\bbX$ (over $\Zp$). The action above gives rise to a homomorphism of Lie algebras

\begin{numequation}\label{map_Lie}
\frg_\Zp \ra H^0(\bbX,\cT_\bbX) \;,
\end{numequation}

\vskip8pt

which is explicitly given by

$$\begin{array}{rcl}
e & \mapsto & \partial_x \\
h_1 & \mapsto & -x\partial_x \\
h_2 & \mapsto & x\partial_x \\
f & \mapsto & \partial_y \end{array}$$

\vskip8pt

On $\bbX$ we consider the sheaf of differential operators $\cD_\bbX^{(m)}$ as defined in \cite{BerthelotDI}, \cite{Huyghe97}. Sections are locally given as finite sums

$$\sum_\nu \gamma_\nu \frac{q^{(m)}_\nu!}{\nu!} \partial_x^\nu \hskip10pt \mbox{ or } \hskip10pt \sum_\nu \gamma'_\nu \frac{q^{(m)}_\nu!}{\nu!} \partial_y^\nu$$

\vskip8pt

with $\gamma_\nu \in \Zp[x]$ and $\gamma_\nu' \in \Zp[y]$, respectively. The sheaf $\cD^{(m)}_\bbX$ is filtered by subsheaves $\cD^{(m)}_{\bbX,d}$ of differential operators of degree $\le d$. Furthermore, for the symmetric algebra $\Sym(\cT_\bbX) = \bigoplus_{d \ge 0} \cT_\bbX^{\otimes d}$ there exists a divided power version

$$\Sym(\cT_\bbX)^{(m)} = \bigoplus _{d \ge 0} (\cT_\bbX^{\otimes d})^{(m)} \;,$$

\vskip8pt

cf. \cite{Huyghe97}. The sheaf $(\cT_\bbX^{\otimes d})^{(m)}$ in degree $d$ is, as $\cO_\bbX$-module, locally generated by

\begin{numequation}\label{local_generators}
\frac{q^{(m)}_{i_1}!}{i_1!} s_1^{\otimes i_1} \cdot \ldots \cdot \frac{q^{(m)}_{i_r}!}{i_r!} s_r^{\otimes i_r} \;,
\end{numequation}

\vskip8pt

where $i_1 + \ldots + i_r = d$ and $s_1, \ldots, s_r$ are local sections of $\cT_\bbX$. There is an obvious monomorphism of sheaves

\begin{numequation}\label{mono_of_sheaves}
\Sym(\cT_\bbX)^{(m)} \hra \Sym(\cT_\bbX)^{(0)}_\bbQ = \Sym(\cT_\bbX)^{(0)} \otimes_\bbZ \bbQ \;.
\end{numequation}

\vskip8pt

\begin{lemma}\label{trivialization} The image of the subsheaf

$$(\cT_\bbX^{\otimes d})^{(m)} \sub \Sym(\cT_\bbX)^{(m)}$$

\vskip8pt

under the map \ref{mono_of_sheaves} is equal to

$$\frac{q^{(m)}_d!}{d!}\cT_\bbX^{\otimes d} \sub
\Sym(\cT_\bbX)^{(0)}_\bbQ \;.$$

\vskip8pt

Therefore,

$$\Sym(\cT_\bbX)^{(m)} = \bigoplus _{d \ge 0} \frac{q^{(m)}_d!}{d!}\cT_\bbX^{\otimes d} \;.$$

\vskip8pt
\end{lemma}

\Pf Because $\cT_\bbX^{\otimes d}$ is locally free of rank one, we can write the local sections $s_i$ in \ref{local_generators} as $s_i = f_i \cdot s$ with a local generator $s$ of $\cT^{\otimes d}$ and local sections $f_i$ of $\cO_\bbX$. Hence we assume $s_i = s$ for $i= 1, \ldots, r$. Moreover, for any $i, j \ge 0$ one has that

\begin{numequation}\label{integrality}
\frac{(i+j)!}{i!j!}
\left(\frac{q^{(m)}_{i+j}!}{q^{(m)}_{i}!q^{(m)}_{j}!}\right)^{-1} \in \Zp \;,
\end{numequation}

cf. \cite[sec. 1]{Huyghe97}. Applying this fact repeatedly shows that

$$\frac{q^{(m)}_{i_1}!}{i_1!} \cdot \ldots \cdot \frac{q^{(m)}_{i_r}!}{i_r!}  \in \frac{q^{(m)}_d!}{d!} \Zp \;,$$

\vskip8pt

and this proves the assertion of the lemma. \qed

\begin{lemma} Fix $d \ge 1$. The map sending $\frac{q^{(m)}_d!}{d!}\partial_x^d$ (resp. $\frac{q^{(m)}_d!}{d!}\partial_y^d$), considered as a local generator of $\cD^{(m)}_{\bbX,d}$ to $\frac{q^{(m)}_d!}{d!}\partial_x^{\otimes d}$ (resp. $\frac{q^{(m)}_d!}{d!}\partial_y ^{\otimes d}$), considered as a local generator of $(\cT_\bbX^{\otimes d})^{(m)}$, induces a canonical exact sequence of sheaves

\begin{numequation}\label{Huyghe_ex_seq}
0 \ra \cD^{(m)}_{\bbX,d-1} \ra \cD^{(m)}_{\bbX,d} \ra (\cT_\bbX^{\otimes d})^{(m)} \ra 0 \;.
\end{numequation}

\end{lemma}

\Pf This is \cite[1.3.7.3]{Huyghe97}. In the case considered here, it is also an immediate consequence of \ref{trivialization}. \qed

\vskip8pt

\begin{prop}\label{global_sections_n_eq_zero} (a) For all $d \ge 0$ one has $H^1(\bbX,\cD_{\bbX,d}^{(m)}) = 0$.

\vskip8pt

(b) For all $d \ge 1$ the sequence

\begin{numequation}\label{ex_seq_global}
0 \ra H^0\left(\bbX,\cD^{(m)}_{\bbX,d-1}\right) \ra H^0\left(\bbX,\cD^{(m)}_{\bbX,d}\right) \ra H^0\left(\bbX,(\cT_\bbX^{\otimes d})^{(m)}\right) \ra 0
\end{numequation}

induced by \ref{Huyghe_ex_seq} is exact.

\vskip8pt

(c) The canonical map

$${\rm gr}\left(H^0\Big(\bbX,\cD_\bbX^{(m)}\Big)\right) = \bigoplus_{d \ge 0} H^0\left(\bbX,\cD^{(m)}_{\bbX,d}\right)\Big/H^0\left(\bbX,\cD^{(m)}_{\bbX,d-1}\right) \;\; \lra \;\; H^0\Big(\bbX,\Sym(\cT_\bbX)^{(m)}\Big)$$

\vskip8pt

is an isomorphism.
\end{prop}

\Pf (a) The proof proceeds by induction on $d$. We have $\cD^{(m)}_{\bbX,0} = \cO_\bbX$, and the assertion is true for $d=0$. Moreover, $\cT_\bbX^{\otimes d} \simeq \cO_\bbX(2d)$ and therefore $H^1(\bbX, \cT_\bbX^{\otimes d}) = 0$. Using \ref{trivialization}, we find that $H^1(\bbX, (\cT_\bbX^{\otimes d})^{(m)}) = 0$ for all $d,m \ge 0$. Now suppose $d \ge 1$. By \ref{Huyghe_ex_seq} we get an exact sequence

$$H^1\left(\bbX,\cD^{(m)}_{\bbX,d-1}\right) \ra H^1\left(\bbX,\cD^{(m)}_{\bbX,d}\right) \ra H^1\left(\bbX,(\cT_\bbX^{\otimes d})^{(m)}\right) \;,$$

\vskip8pt

and our induction hypothesis implies $H^1\left(\bbX,\cD^{(m)}_{\bbX,d}\right) = 0$.

\vskip8pt

(b) This assertion follows from (a) and the long exact cohomology sequence attached to \ref{Huyghe_ex_seq}.

\vskip8pt

(c) This follows immediately from (b). \qed

\begin{rem} Assertion (c) of the previous proposition is as in \cite[2.3.6 (ii)]{Huyghe97}, at least for large $d$. Though Noot-Huyghe's result would be good enough for our purposes, we have preferred to give a self-contained proof here. The proof given here proceeds along the same lines as the proof in \cite{Huyghe97}. \qed
\end{rem}

\vskip8pt

In the following we consider the filtration of $U(\frg_\Zp)^{(m)}$ whose submodule of elements of degree $\le d$ is generated as a $\Zp$-module by terms of the form \ref{generators} with $\nu_1+\nu_2+\nu_3+\nu_4 \le d$.

\vskip8pt

\begin{prop}\label{map} (a) For all $\nu \ge 0$ one has the following identity of differential operators in $\cD_\bbX \otimes_\bbZ \bbQ$: ${x\partial_x \choose \nu} = x^\nu \frac{\partial_x^\nu}{\nu!}$.

\vskip8pt

(b) The canonical map $U(\frg_\Zp) \ra H^0(\bbX,\cD^{(0)}_\bbX)$ induced by \ref{map_Lie} extends to a homomorphism

\begin{numequation}\label{map_U}
\xi^{(m)}: U(\frg_\Zp)^{(m)} \lra H^0(\bbX,\cD^{(m)}_\bbX) \;,
\end{numequation}

of $\Zp$-algebras which is compatible with the filtrations on both sides.

\vskip8pt

(c) $\xi^{(m)}$ maps the center $Z(\frg_\Zp)$ of $U(\frg_\Zp) \sub U(\frg_\Zp)^{(m)}$ to $\Zp$. Let $\theta_0 = \xi^{(m)}|_{Z(\frg_\Zp)}$ be the restriction of $\xi^{(m)}$ to $Z(\frg_\Zp)$. Then $\ker(\xi^{(m)})$ is the (two-sided) ideal generated by $\ker(\theta_0)$.

\end{prop}

\Pf (a) Is easily proved by induction.

\vskip8pt

(b) Using (a) we see that ${h_i \choose \nu}$, $i=1,2$, is mapped to $\pm x^\nu \frac{\partial_x^\nu}{\nu!}$. The assertion now follows directly from the definition of $U(\frg_\Zp)^{(m)}$.

\vskip8pt

(c) Tensor with $\bbQ$ and use the statement in characteristic zero, cf. \cite{BB81}. \qed

\vskip8pt

Using the notations introduced in \ref{map} we define

$$U(\frg_\Zp)^{(m)}_{\theta_0} \stackrel[\rm df]{}{=} U(\frg_\Zp)^{(m)} \otimes_{Z(\frg_\Zp),\theta_0} \Zp \;.$$

\vskip8pt

Therefore, $\xi^{(m)}$ induces an injective homomorphism of $\Zp$-algebras

\begin{numequation}\label{map_U_0}
\xi^{(m)}_0: U(\frg_\Zp)^{(m)}_{\theta_0} \hra H^0(\bbX,\cD^{(m)}_\bbX) \;.
\end{numequation}

\begin{prop}\label{cokernel_bdd_torsion} (a) Via the homomorphism

$${\rm gr}\, \xi^{(m)}: {\rm gr}\Big(U(\frg_\Zp)^{(m)}\Big) \lra H^0\Big(\bbX,{\rm Sym}(\cT_\bbX)^{(m)}\Big) = {\rm gr}\left(H^0\Big(\bbX,\cD_\bbX^{(m)}\Big)\right)$$

\vskip8pt

induced by $\xi^{(m)}$, the ring $H^0\Big(\bbX,{\rm Sym}(\cT_\bbX)^{(m)}\Big)$ is a finitely generated module over \linebreak ${\rm gr}\Big(U(\frg_\Zp)^{(m)}\Big)$.

(b) Via $\xi^{(m)}_0$ the ring $H^0(\bbX,\cD^{(m)}_\bbX)$ is a finitely generated $U(\frg_\Zp)^{(m)}_{\theta_0}$-module. Moreover, there is $N(m) \in \bbN$ such that the $\coker(\xi^{(m)}_0)$ is annihilated by $p^{N(m)}$.
\end{prop}

\Pf (a) By \ref{trivialization} we have

$$H^0\left(\bbX,\Sym(\cT_\bbX)^{(m)}\right) = \bigoplus _{d \ge 0} \frac{q^{(m)}_d!}{d!}H^0\left(\bbX,\cT_\bbX^{\otimes d}\right) \;,$$

\vskip8pt

as submodules of $H^0\left(\bbX,\Sym(\cT_\bbX)^{(m)}_\bbQ \right)= H^0\left(\bbX,\Sym(\cT_\bbX)^{(m)}\right) \otimes_\bbZ \bbQ$.
Furthermore,

$$H^0\left(\bbX,\cT_\bbX^{\otimes d}\right) = \bigoplus_{k=0}^{2d} \Zp x^k \partial_x^{\otimes d} \;.$$

\vskip8pt

Our goal is to show that $H^0\left(\bbX,\Sym(\cT_\bbX)^{(m)}\right)$ is generated as a module over ${\rm gr}\Big(U(\frg_\Zp)^{(m)}\Big)$ by the elements

$$\frac{q^{(m)}_d!}{d!} x^k \partial_x^{\otimes d} \hskip10pt \mbox{ with
} \hskip4pt 0 \le d < 2p^m \;, \;\; 0 \le k \le 2d \;.$$

\vskip8pt

To this end, consider an element $\frac{q^{(m)}_d!}{d!} x^k \partial_x^{\otimes d}$ with $k \le 2d$. Write $d =  p^mq + s$. We are going to use the elementary fact

$$\frac{q^{(m)}_d!}{d!} = \frac{u}{s!(p^m!)^q} \;,$$

with a $p$-adic unit $u$, cf. \cite[5.2.2]{EmertonA}.

\vskip8pt

{\it Case $k \le d$.} Writing $k = p^mq' + r$, we have $q' \le q$. If $r \le 2s$ then consider the equation\footnote{This equation and the following formulas are to be considered in the commutative ring $H^0\left(\bbX,\Sym(\cT_\bbX)^{(m)}\right)$. To simplify notation we have dropped the superscript ''$\otimes$''.}

$$\frac{q^{(m)}_d!}{d!} x^k \partial_x^d = u\left(\frac{(x\partial_x)^{p^m}}{p^m!}\right)^{q'} \cdot \left(\frac{\partial_x^{p^m}}{p^m!}\right)^{q-q'} \cdot \frac{1}{s!}x^r\partial_x^s \;.$$

\vskip8pt

Now suppose $r>2s$. Because $k=p^mq'+r \le d = p^mq+s$ we must have $q'<q$ and hence $q-q'-1\ge 0$. Then we can write

$$\frac{q^{(m)}_d!}{d!} x^k \partial_x^d = u\left(\frac{(x\partial_x)^{p^m}}{p^m!}\right)^{q'} \cdot \left(\frac{\partial_x^{p^m}}{p^m!}\right)^{q-q'-1} \cdot \frac{1}{s!(p^m!)}x^r\partial_x^{p^m+s} \;.$$

\vskip8pt

{\it Case $d < k$ ($\le 2d$).} Write $k=p^mq'+r$, and suppose $q'=2q''$ is even. Because $\frac{k}{2} = p^mq''+\frac{r}{2} \le p^mq+s$ we must have $q'' \le q$. If $r \le 2s$ then consider the equation

$$\frac{q^{(m)}_d!}{d!} x^k \partial_x^d = u\left(\frac{(x^2\partial_x)^{p^m}}{p^m!}\right)^{q''} \cdot \left(\frac{\partial_x^{p^m}}{p^m!}\right)^{q-q''} \cdot \frac{1}{s!}x^r\partial_x^s \;.$$

\vskip8pt

Now suppose $r > 2s$. Then we must have $q''<q$, hence $q-q''-1 \ge 0$ and we can write

$$\frac{q^{(m)}_d!}{d!} x^k \partial_x^d = u\left(\frac{(x^2\partial_x)^{p^m}}{p^m!}\right)^{q''} \cdot \left(\frac{\partial_x^{p^m}}{p^m!}\right)^{q-q''-1} \cdot \frac{1}{s!(p^m!)}x^{r}\partial_x^{p^m+s} \;.$$

\vskip8pt

Assume now that $q'=2q''+1$ is odd. Because

$$p^mq''+\frac{p^m+r}{2} = p^m(q''+\frac{1}{2})+\frac{r}{2}  = \frac{k}{2} \le d = p^mq+s \;,$$

\vskip8pt

we must have $q'' \le q$. If $p^m+r \le 2s$ we consider

$$\frac{q^{(m)}_d!}{d!} x^k \partial_x^d = u\left(\frac{(x^2\partial_x)^{p^m}}{p^m!}\right)^{q''} \cdot \left(\frac{\partial_x^{p^m}}{p^m!}\right)^{q-q''} \cdot\frac{1}{s!}x^{p^m+r}\partial_x^s \;.$$

\vskip8pt

Finally, if $p^m+r > 2s$ we must have $q''<q$. In this case we consider

$$\frac{q^{(m)}_d!}{d!} x^k \partial_x^d = u\left(\frac{(x^2\partial_x)^{p^m}}{p^m!}\right)^{q''} \cdot \left(\frac{\partial_x^{p^m}}{p^m!}\right)^{q-q''-1} \cdot\frac{1}{s!(p^m!)}x^{p^m+r}\partial_x^{p^m+s} \;.$$

\vskip8pt

(b) For $0 \le d < 2p^m$ and $0 \le k \le 2d$ let $e_{d,k}$ be a representative in $H^0\left(\bbX,\cD_{\bbX,d}^{(m)}\right)$ of the element $\frac{q^{(m)}_d!}{d!} x^k \partial_x^{\otimes d}$ in $H^0\left(\bbX,(\cT_\bbX^{\otimes d})^{(m)}\right)$.
By part (a), $H^0\Big(\bbX,{\rm Sym}(\cT_\bbX)^{(m)}\Big)$ is generated over ${\rm gr}\Big(U(\frg_\Zp)^{(m)}\Big)$ by the elements $\frac{q^{(m)}_d!}{d!} x^k \partial_x^{\otimes d}$, for $0 \le d < 2p^m$ and $0 \le k \le 2d$, it follows that
$H^0(\bbX,\cD^{(m)}_\bbX)$ is generated over $U(\frg_\Zp)^{(m)}$ by the elements $e_{d,k}$. And then, obviously, $H^0(\bbX,\cD^{(m)}_\bbX)$ is actually a finitely generated $U(\frg_\Zp)^{(m)}_{\theta_0}$-module. Moreover, we see that the generators

$$\frac{q^{(m)}_d!}{d!} x^k \partial_x^{\otimes d} \hskip10pt \mbox{ with
} \hskip4pt 0 \le d < 2p^m \;, \;\; 0 \le k \le 2d \;, $$

\vskip8pt

of $H^0\Big(\bbX,{\rm Sym}(\cT_\bbX)^{(m)}\Big)$ over ${\rm gr}\Big(U(\frg_\Zp)^{(m)}\Big)$ have the property that

\begin{numequation}\label{torsion}
(p^m-1)! \cdot (p^m)! \cdot \frac{q^{(m)}_d!}{d!} x^k \partial_x^{\otimes d} \in \im\left({\rm gr}\Big(U(\frg_\Zp)^{(m)}\Big) \ra H^0\Big(\bbX,{\rm Sym}(\cT_\bbX)^{(m)}\Big)\right) \;.
\end{numequation}

Because the generators $e_{d,k}$ are in degrees $< 2p^m$, repeating \ref{torsion} finitely often shows that there is $N(m)$ such that $p^{N(m)}e_{d,k}$ is in the image of $U(\frg_\Zp)^{(m)}_{\theta_0}$ for $0 \le d < 2p^m$, $0 \le k \le 2d$. Now assertion (b) follows. \qed

\vskip8pt

\subsection{$\sD^\dagger$ and the distribution algebra $\cD^{\rm an}(\bbG(0)^\circ)$}

Denote by $\frX$ the completion of $\bbX$ along its special fiber $\bbX_\Fp$. Let $\sD_\frX^{(m)}$ be the $p$-adic completion of the sheaf $\cD_\bbX^{(m)}$, which we consider as a sheaf on $\frX$.

\begin{lemma}\label{completion} The canonical map

$$H^0\left(\bbX,\cD^{(m)}_\bbX\right) \lra H^0\left(\frX,\sD^{(m)}_\frX\right)$$

\vskip8pt

extends to an isomorphism

$$H^0\left(\bbX,\cD^{(m)}_\bbX\right)^\wedge \lra H^0\left(\frX,\sD^{(m)}_\frX\right) \;,$$

\vskip8pt

where the left hand side is the $p$-adic completion of $H^0\left(\bbX,\cD^{(m)}_\bbX\right)$.

\end{lemma}

\Pf This is contained in \cite[Prop. 3.2]{Huyghe97}. The key ingredient used in \cite[Prop. 3.2]{Huyghe97} is that $H^1$ of the sheaf in question (here $\cD^{(m)}_\bbX$) is annihilated by a finite power of $p$. Here we have seen $H^1(\bbX,\cD_\bbX^{(m)}) = 0$, cf. \ref{global_sections_n_eq_zero}. Thus it would be possible to give a self-contained proof following the arguments given in the proof of \cite[Prop. 3.2]{Huyghe97}. \qed

\vskip8pt

Put

$$\sD^\dagger_\frX = \varinjlim_m \sD_\frX^{(m)} \;,$$

\vskip8pt

and

$$\sD^\dagger_{\frX,\bbQ} = \varinjlim_m \sD_\frX^{(m)}\otimes_\bbZ \bbQ \;.$$

\vskip8pt

\begin{thm}\label{thm-0} (a) The homomorphism

$$\xi^{(m)}_0: U\left(\frg_\Zp\right)^{(m)}_{\theta_0} \ra H^0\left(\bbX,\cD^{(m)}_\bbX\right) \;,$$

\vskip8pt

cf. \ref{map_U_0}, induces a homomorphism

$$\widehat{\xi}^{(m)}_0: \widehat{U}\left(\frg_\Zp\right)^{(m)}_{\theta_0} \ra H^0\left(\frX,\sD^{(m)}_\frX\right) \;,$$

\vskip8pt

which is injective and whose cokernel is annihilated by $p^{N(m)}$ where $N(m)$ is as in \ref{cokernel_bdd_torsion}. Therefore, $\widehat{\xi}^{(m)}_0$ induces an isomorphism

$$\widehat{U}\left(\frg_\Zp\right)^{(m)}_{\theta_0,\bbQ} \stackrel{\simeq}{\lra} H^0\left(\frX,\sD^{(m)}_{\frX,\bbQ}\right) \;.$$

\vskip8pt

(b) The isomorphisms in (a) give rise to a canonical isomorphism

$$ \cD^{\rm an}(\bbG(0)^\circ)_{\theta_0} \simeq \varinjlim_m \widehat{U}\left(\frg_\Zp\right)^{(m)}_{\theta_0,\bbQ} \stackrel{\simeq}{\lra} H^0(\frX,\sD^\dagger_{\frX,\bbQ}) \;.$$
\end{thm}

\Pf (a) We consider the exact sequence induced by $\xi^{(m)}_0$

\vskip12pt

$0 \lra \left(U\left(\frg_\Zp\right)^{(m)}_{\theta_0} \cap p^kH^0\left(\bbX,\cD^{(m)}_\bbX\right)\right) \Big/ p^kU\left(\frg_\Zp\right)^{(m)}_{\theta_0}$

\vskip8pt

\hfill $\lra U\left(\frg_\Zp\right)^{(m)}_{\theta_0} \Big/ p^kU\left(\frg_\Zp\right)^{(m)}_{\theta_0} \lra H^0\left(\bbX,\cD^{(m)}_\bbX\right) \Big/ p^kH^0\left(\bbX,\cD^{(m)}_\bbX\right) \;$.

\vskip8pt

Because the projective limit functor is left-exact, and as $H^0\left(\bbX,\cD^{(m)}_\bbX\right)$ is separated for the $p$-adic topology, we deduce that the homomorphism $\widehat{\xi}^{(m)}_0$ between the completions is injective as well. The assertion about the cokernel is an immediate consequence of \ref{cokernel_bdd_torsion}. Hence the isomorphism after extending scalars to $\bbQ$.

\vskip8pt

(b) This assertion follows from (a) and the fact that cohomology commutes with direct limits. \qed

\vskip8pt

As already indicated in the introduction, after having obtained
this result we have been informed by C. Noot-Huyghe that she has
proved the general case of this theorem, for an arbitrary split
reductive group and the corresponding flag variety, in an
unpublished manuscript.

The isomorphism in (a) for an arbitrary split semisimple group and
the corresponding flag variety has appeared, in the case $m=0$ and
with some restrictions on the prime number $p$, in \cite{AW}.

\section{The semistable models $\bbX_n$ and their completions $\frX_n$}\label{Xn}

\subsection{The construction via blowing-up}\label{blowing_up}

\begin{para}\label{construction_of_X_n} In the following, all closed subsets of a scheme are considered as closed subschemes with their reduced induced subscheme structure. Put $\bbX_0 = \bbX = \bbP^1_\Zp$. Blowing up $\bbX_0$ in the $\Fp$-rational points of its special fiber $\bbX_{0,\Fp}$ produces a scheme $\bbX_1$. The irreducible components of the special fiber of $\bbX_1$ are all projective lines over $\Fp$, and there are $p+2$ of them: on the one hand we have the strict transform of $\bbX_{0,\Fp}$, which we can and will identify with $\bbX_{0,\Fp}$, and then there is for any $\Fp$-rational point $P$ of $\bbX_0$ the corresponding component $E_P \simeq \bbP^1_\Fp$ of the exceptional divisor. No two components $E_P$ intersect each other, but any one of these intersects $\bbX_{0,\Fp}$ in a unique point
which corresponds to the point $P$ that has been blown up. We call the components $E_P$ the 'end components' or 'ends' of the special fiber of $\bbX_1$.

\vskip8pt

Then blow up $\bbX_1$ in the smooth $\Fp$-rational points of its special fiber. There are $p$ such points on each component $E_P$. Call the resulting scheme $\bbX_2$. The special fiber of $\bbX_2$ consists of the strict transform of the special fiber of $\bbX_1$, which we identify with $\bbX_{1,\Fp}$, and, for each of the components $E_P$ of $\bbX_{1,\Fp}$ there are $p$ irreducible components $E_{P,P'}$ of the exceptional divisor, and $E_{P.P'}$ intersects $E_P$ in the point $P'$ that has been blown up. Again, we call the irreducible components $E_{P,P'}$ the 'end components' or 'ends' of the special fiber of $\bbX_2$.

\vskip8pt

Inductively one defines $\bbX_n$ by blowing up $\bbX_{n-1}$ in the
smooth $\Fp$-rational points of the special fiber of $\bbX_{n-1}$.
The irreducible components of the exceptional divisor are called
the 'end components' or 'ends' of the special fiber of $\bbX_n$.
It is easy to see that the intersection graph of the special fiber
of $\bbX_n$ is a tree. There are $p+1$ edges meeting at every
vertex, except for the vertices which correspond to the end
components: these are only connected to the rest of the tree by a
single edge.
\end{para}

\vskip8pt

\begin{rem}\hskip-4pt\footnote{The content of this remark will not be used later on.} Because the group $\bbG(\Zp) = \GL_2(\Zp)$ acts on $\bbX_0$ and preserves the closed subscheme $\bbX_0(\Fp)$, the group $\bbG(\Zp)$ acts also on $\bbX_1$. It is easy to see that $\bbG(\Zp)$ preserves the subscheme of $\bbX_1$ which gives rise to $\bbX_2$. Inductively we find that $\bbG(\Zp)$ acts on $\bbX_n$ for all $n$. Furthermore, one can show that the group scheme $\bbG(n)$ acts on the scheme $\bbX_n$.
\end{rem}

\vskip8pt

\subsection{An open affine covering of $\bbX_n$}\label{covering}

Here we describe an open affine covering of the scheme $\bbX_n$, and a coherent system of local coordinates\footnote{By this we mean a set of local coordinates together with transition formulas for the local coordinates on 'neighboring' open affine subsets. The meaining of 'neighboring' in our context will become clear in the sequel.}. This will be used later in sec. \ref{log_tangent_sheaf}.

\begin{para}\label{general_shape}{\it Outline.} We will first describe the general shape of this covering and the procedure by which it is obtained. Let $\cR \sub \Zp$ be any system of representatives for $\Zp/p\Zp$ and put $\cR_\infty = \cR \cup \{\infty\}$. Let $n \ge 1$. Inductively we will
define an open subset $\bbX_{n-1}^\circ \sub \bbX_{n-1}$ and open affine 'residual disc schemes' $\bbD^{(n-1)}_{\underline{a}}$ for any tuple $\underline{a} = (a_0,a_1, \ldots, a_{n-1}) \in \cR_\infty \times \cR^{n-1}$. Each scheme $\bbD^{(n-1)}_{\underline{a}}$ has a unique $\Fp$-rational point and $\bbX_n$ is obtained from $\bbX_{n-1}$ by blowing up all these points. The open subset
$\bbX_{n-1}^\circ \sub \bbX_{n-1}$ is not affine (except if $n-1=0$) but it is equipped with an open affine covering. Moreover, the special fiber of $\bbX_{n-1}^\circ$ does not contain any smooth $\Fp$-rational point of the special fiber of $\bbX_{n-1}$. The blow-up morphism $pr_{n,n-1}: \bbX_n \ra \bbX_{n-1}$ is thus an isomorphism over $\bbX_{n-1}^\circ$, and the preimage $pr_{n,n-1}^{-1}(\bbX_{n-1}^\circ) \sub \bbX_n$ is then equipped with the open affine covering of $\bbX_{n-1}^\circ$. In the following we identify $pr_{n,n-1}^{-1}(\bbX_{n-1}^\circ)$ with $\bbX_{n-1}^\circ$.

\vskip8pt

Next we define for any such $\underline{a}$ open affine subschemes $\bbX^{(n)}_{\underline{a}}$ and, for all $a_n \in \cR$, 'residual disc schemes' $\bbD^{(n)}_{\underline{a},a_n}$ of $\bbX_n$. These open affine subschemes, together with the open affine covering of $\bbX_{n-1}^\circ$ constitute then the open affine covering of $\bbX_n$. The open subset $\bbX_n^\circ$ is defined as

$$\bbX_n^\circ  \stackrel[\rm df]{}{=} \bbX_{n-1}^\circ \cup \bigcup_{\underline{a} \in \cR_\infty \times \cR^{n-1}}\bbX^{(n)}_{\underline{a}} \;.$$

\vskip8pt
\end{para}

\begin{para}\label{bbX_0}{\it When $n=0$.} We start with the affine covering $\bbX_0 = U_x \cup U_y$ of $\bbX_0$, cf. \ref{smooth}, where
$U_x = \Spec(\Zp[x])$ and $U_y = \Spec(\Zp[y])$ and these open subschemes are glued together according to the relation $xy=1$. For $a \in \cR$ put $x^{(0)}_a = x-a$, and consider this as a local coordinate at $x = a$, and set $x^{(0)}_\infty = y = \frac{1}{x}$. For $a \in \cR_\infty$ put

$$R^{(0)}_a = \Zp[x^{(0)}_a]\left[\frac{1}{x^{(0)}_b} \;\; \Bigg| \;\; b \in \cR, b \neq a \right] \;,$$

\vskip8pt

and view this as a subring of the rational function field $\Qp(x)$. It is immediate that for all $a \in \cR_\infty$ the ring

$$R^{(0)} \stackrel[\rm df]{}{=} R^{(0)}_a\left[\frac{1}{x^{(0)}_a}\right] \;,$$

\vskip8pt

as a subring of $\Qp(x)$, is independent of $a$. Set $\bbX_0^\circ = \Spec(R^{(0)})$. The special fiber of $\bbX_0^\circ$ is $\bbP^1_\Fp \setminus \bbP^1(\Fp)$. Furthermore, for $a \in \cR_\infty$ put

$$\bbD^{(0)}_a = \Spec(R^{(0)}_a)  \;.$$

\vskip8pt

The special fiber of $\bbD^{(0)}_a$ is $\bbX_{0,\Fp}^\circ \cup \{\overline{a}\}$, where $\overline{a}$ is the 'mod-$p$ reduction of $a$'. This is  the unique $\Fp$-rational point which corresponds to the ideal $(p,x^{(0)}_a)$. We call $\bbD^{(0)}_a$ a 'residual disc scheme'. For later use we fix the coordinate function $x^{(0)}_a$ on $\bbD^{(0)}_a$. For any two distinct $a,a' \in \cR_\infty$ we have $\bbD^{(0)}_a \cap \bbD^{(0)}_{a'} = \bbX_0^\circ$. Then we consider the covering of $\bbX_0$ by the open subschemes $\bbD^{(0)}_a$, $a \in \cR_\infty$, together with $\bbX_0^\circ$.
\end{para}

\vskip8pt

\begin{para}\label{bbX_1}{\it When $n=1$.} $\bbX_1$ is obtained by blowing up $\bbX_0$ in the points corresponding to the ideals $(p,x^{(0)}_{a_0}) \sub R^{(0)}_{a_0}$, $a_0 \in \cR_\infty$. In order to describe $\bbX_1$, we introduce new indeterminates $z^{(1)}_{a_0}$ and $x^{(1)}_{a_0}$ satisfying

$$x^{(0)}_{a_0}z^{(1)}_{a_0} = p \hskip12pt \mbox{and} \hskip12pt z^{(1)}_{a_0} x^{(1)}_{a_0} = 1 \;.$$

\vskip8pt

Set also $x^{(1)}_{a_0,a_1} = x^{(1)}_{a_0}- a_1$ for $a_1 \in \cR$. Then define

$$R^{(1)}_{a_0} = R^{(0)}_{a_0}[z^{(1)}_{a_0}]\left[\frac{1}{x^{(1)}_{a_0,a_1}} \;\; \Bigg| \;\; a_1 \in \cR \right] \Bigg/(x^{(0)}_{a_0}z^{(1)}_{a_0}-p)\;,$$

\vskip8pt

and put $\bbX^{(1)}_{a_0} = \Spec(R^{(1)}_{a_0})$. For $a_1 \in \cR$ set

$$R^{(1)}_{a_0,a_1} =  R^{(0)}_{a_0}[x^{(1)}_{a_0,a_1}]\left[\frac{1}{x^{(1)}_{a_0,b}} \;\; \Bigg| \;\; b \in \cR \setminus \{a_1\}\right] \;,$$

and define

$$\bbD^{(1)}_{a_0,a_1} = \Spec\left(R^{(1)}_{a_0,a_1}\right) \;.$$

\vskip8pt

The special fiber of each $\bbD^{(1)}_{a_0,a_1}$ is isomorphic to an affine line over $\Fp$ all of whose $\Fp$-rational points have been removed, except one. Again, in order to obtain a coherent system of coordinates, we fix the coordinate function $x^{(1)}_{a_0,a_1}$ on $\bbD^{(1)}_{a_0,a_1}$. For any $a_1 \in \cR$ one has

$$R^{(1)}_{a_0}\left[\frac{1}{z^{(1)}_{a_0}}\right] = R^{(1)}_{a_0,a_1}\left[\frac{1}{x^{(1)}_{a_0,a_1}}\right] \;,$$

\vskip8pt

and this ring is thus independent of $a_1$. For any two distinct $a_1,a_1' \in \cR$ one has

$$\bbD^{(1)}_{a_0,a_1} \cap \bbD^{(1)}_{a_0,a_1'} = \bbD^{(1)}_{a_0,a_1} \cap \bbX^{(1)}_{a_0}  = \Spec\left(R^{(1)}_{a_0}\left[\frac{1}{z^{(1)}_{a_0}}\right]\right)\;,$$

\vskip8pt

and the special fiber of this scheme is isomorphic (via the coordinate $x^{(1)}_{a_0}$, say) to $\bbP^1_\Fp \setminus \bbP^1(\Fp)$.
Furthermore, for any two distinct $a_0,a_0' \in \cR_\infty$ one has

$$\bbX^{(1)}_{a_0} \cap \bbX^{(1)}_{a_0'} = \bbX_0^\circ \;.$$

\vskip8pt

Let $\bbX_1^\circ$ be the union of the schemes $\bbX^{(1)}_{a_0}$, $a_0 \in \cR_\infty$, and $\bbX_0^\circ$.
\end{para}

\begin{para}\label{bbX_n}{\it From $n-1$ to $n$.}
Firstly, we use the preimages of the affine covering of $\bbX_{n-1}^\circ$ under the blow-up map $\bbX_n \ra \bbX_{n-1}$. Then we consider a 'residue disc scheme'

$$\bbD^{(n-1)}_{\underline{a}} = \Spec\left(R^{(n-1)}_{\underline{a}}\right)$$

\vskip8pt

of $\bbX_{n-1}$, where $\underline{a} = (a_0,a_1,\ldots,a_{n-1})$. It is equipped with a coordinate function $x^{(n-1)}_{\underline{a}}$ and has a unique $\Fp$-rational point which corresponds to the ideal $(p,x^{(n-1)}_{\underline{a}}) \sub R^{(n-1)}_{\underline{a}}$.
$\bbX_n$ is obtained from $\bbX_{n-1}$ by blowing up these $\Fp$-rational points, for all $\underline{a} \in \cR_\infty \times \cR^{n-1}$.

\vskip8pt

To describe the blow-up process, we introduce indeterminates $z^{(n)}_{\underline{a}}$ and $x^{(n)}_{\underline{a}}$ satisfying

$$x^{(n-1)}_{\underline{a}}z^{(n)}_{\underline{a}} = p \hskip12pt \mbox{and} \hskip12pt z^{(n)}_{\underline{a}} x^{(n)}_{\underline{a}} = 1 \;.$$

\vskip8pt

For $a_n \in \cR$ set $x^{(n)}_{\underline{a},a_n} = x^{(n)}_{\underline{a}}-a_n$ and define

\begin{numequation}\label{local_eqn} R^{(n)}_{\underline{a}} = R^{(n-1)}_{\underline{a}}[z^{(n)}_{\underline{a}}]\left[\frac{1}{x^{(n)}_{\underline{a},b}} \;\; \Bigg| \;\; b \in \cR \right] \Bigg/ (x^{(n-1)}_{\underline{a}}z^{(n)}_{\underline{a}} - p) \;,
\end{numequation}

and put

$$\bbX^{(n)}_{\underline{a}} = \Spec\left(R^{(n)}_{\underline{a}}\right) \;.$$

\vskip8pt

For $a_n \in \cR$ define

$$R^{(n)}_{\underline{a},a_n} = R^{(n-1)}_{\underline{a}}[x^{(n)}_{\underline{a},a_n}]\left[\frac{1}{x^{(n)}_{\underline{a},b}} \;\; \Bigg| \;\; b \in \cR \setminus \{a_n\} \right] \;,$$

\vskip8pt

and put

$$\bbD^{(n)}_{\underline{a},a_n} = \Spec\left(R^{(n)}_{\underline{a},a_n}\right) \;.$$

\vskip8pt

Again, in order to obtain a coherent system of coordinates, we fix the coordinate function $x^{(n)}_{\underline{a},a_n}$ on $\bbD^{(n)}_{\underline{a},a_n}$. For any $a_n \in \cR$ one has

$$R^{(n)}_{\underline{a}}\left[\frac{1}{z^{(n)}_{\underline{a}}}\right] = R^{(n)}_{\underline{a},a_n}\left[\frac{1}{x^{(n)}_{\underline{a},a_n}}\right] \;,$$

\vskip8pt

and this ring is thus independent of $a_n$. For any two distinct $a_n,a_n' \in \cR$ one has

$$\bbD^{(n)}_{\underline{a},a_n} \cap \bbD^{(n)}_{\underline{a},a_n'} = \bbD^{(n)}_{\underline{a},a_n} \cap \bbX^{(n)}_{\underline{a}}  = \Spec\left(R^{(n)}_{\underline{a}}\left[\frac{1}{z^{(n)}_{\underline{a}}}\right]\right)\;,$$

\vskip8pt

and the special fiber of this scheme is isomorphic to (via the coordinate $x^{(n)}_{\underline{a}}$, say) to $\bbP^1_\Fp \setminus \bbP^1(\Fp)$. Let $\bbX_n^\circ$ be the union of the schemes $\bbX^{(n)}_{\underline{a}}$, $\underline{a} \in \cR_\infty \times \cR^{n-1}$, and $\bbX_{n-1}^\circ$. One obtains an open affine cover for $\bbX_n^\circ$ from the union of the open affine cover from $\bbX_{n-1}^\circ$ and the collection of all $\bbX^{(n)}_{\underline{a}}$. Finally $\bbX_n$ is then covered by $\bbX_n^\circ$ and the open affine subschemes $\bbD^{(n)}_{\underline{a},a_n}$, $(\underline{a},a_n) = (a_0, \ldots, a_{n-1},a_n) \in \cR_\infty \times \cR^n$. Writing out the open affine covering of $\bbX_n^\circ$ explicitly gives:

\begin{numequation}\label{explicit_covering}\bbX_n = \bbX_0^\circ \;\; \cup \; \bigcup_{1 \le \nu \le n} \;\; \bigcup_{\underline{a} \in \cR_\infty \times \cR^{\nu-1}} \bbX^{(\nu)}_{\underline{a}} \;\; \cup \;\; \bigcup_{\underline{b} \in \cR_\infty \times \cR^{n}} \bbD^{(n)}_{\underline{b}} \;.
\end{numequation}

\vskip8pt

\end{para}

\begin{para}\label{formulas} Going through the successive definitions of the local coordinates $x^{(0)}_{a_0}, x^{(1)}_{a_0,a_1}, \ldots, x^{(n)}_{\underline{a}}$, $\underline{a} = (a_0, \ldots, a_n)$, one finds the relations, for $a_0 \neq \infty$,

\begin{numequation}\label{local_x_coord}
\begin{array}{rcl}x & = & a_0+a_1p+ \ldots+a_{n-1}p^{n-1}+ a_np^n + p^nx^{(n)}_{\underline{a}} \\
&&\\
& = & a_0+a_1p+ \ldots+a_{n-1}p^{n-1}+ p^nx^{(n)}_{(a_0, a_1, \ldots, a_{n-1})}\\
&&\\
& = & a_0+a_1p+ \ldots+a_{n-1}p^{n-1}+ p^{n-1}x^{(n-1)}_{(a_0, a_1, \ldots, a_{n-1})} \;,
\end{array}
\end{numequation}

where we have used $x^{(n-1)}_{(a_0, a_1, \ldots, a_{n-1})}z^{(n)}_{(a_0, a_1, \ldots, a_{n-1})} = p$. Similarly we have for \linebreak
$\underline{a} = (\infty, a_1,\ldots, a_{n-1},a_n)$ and $y$ the relations

\begin{numequation}\label{local_y_coord}
\begin{array}{rcl}y & = & a_1p+ \ldots+a_{n-1}p^{n-1}+a_np^n + p^nx^{(n)}_{\underline{a}} \\
&&\\
& = & a_1p+ \ldots+a_{n-1}p^{n-1}+ p^nx^{(n)}_{(a_0, a_1, \ldots, a_{n-1})} \\
&&\\
& = & a_1p+ \ldots+a_{n-1}p^{n-1}+ p^{n-1}x^{(n-1)}_{(a_0, a_1, \ldots, a_{n-1})} \;.
\end{array}
\end{numequation}
\end{para}

\subsection{The formal schemes $\frX_n$}\label{frX_n}

\begin{para} We denote by $\frX_n$ the completion of $\bbX_n$ along its special fiber. One can also obtain $\frX_n$ directly from $\frX$ by the same procedure as in \ref{blowing_up}. Assuming we have constructed $\frX_{n-1}$, we define $\frX_n$ by blowing up (in the sense of formal geometry) the smooth $\Fp$-rational points of the special fiber of $\frX_{n-1}$.

\vskip8pt

Furthermore, the open affine covering described in \ref{covering} gives rise upon completion to a covering of $\frX_n$ by open affine subschemes. The explicit description of the formal completion $\widehat{\bbX}^{(n)}_{\underline{a}}$ of $\bbX^{(n)}_{\underline{a}}$, $\underline{a} \in \cR_\infty \times \cR^{n-1}$ is in fact simpler than the corresponding description for $\bbX^{(n)}_{\underline{a}}$. One can show

$$\widehat{\bbX}^{(n)}_{\underline{a}} = \Spf\left(\Zp\langle x^{(n-1)}_{\underline{a}},z^{(n)}_{\underline{a}}\rangle \left[\frac{1}{(x^{(n-1)}_{\underline{a}})^{p-1}-1}, \frac{1}{(z^{(n)}_{\underline{a}})^{p-1}-1} \right] \Big/(x^{(n-1)}_{\underline{a}}z^{(n)}_{\underline{a}}-p)\right) \;.$$

\vskip8pt

See \cite{TeitelbaumUniversal}\footnote{The relevant material is in the section ''The formal scheme $\widehat{\sH}_p$ -- the naive construction''.} and \cite[I.3]{BoutotCarayol} for details. Similarly, the formal completion $\widehat{\bbD}^{(n)}_{\underline{a},a_n}$ of $\bbD^{(n)}_{\underline{a},a_n}$, $a_n \in \cR$, can be described by

$$\widehat{\bbD}^{(n)}_{\underline{a},a_n} = \Spf\left(\Zp\langle x^{(n)}_{\underline{a},a_n}\rangle \left[\frac{1}{(x^{(n)}_{\underline{a},a_n})^{p-1}-1} \right] \right) \;.$$

\vskip8pt

\end{para}

\begin{rem} Denote by $\frX_n^\circ$ the completion of $\bbX_n^\circ$ along its special fiber. The open embedding $\bbX_{n-1}^\circ \hra \bbX_n^\circ$ induces an open embedding $\frX_{n-1}^\circ \hra \frX_n^\circ$. ($\frX_n^\circ$ can also defined intrinsically, and more straightforwardly, without the use of $\bbX_n^\circ$.) The inductive limit $\varinjlim_n \frX_n^\circ$ is then a formal model of the {\it $p$-adic upper half plane}, cf. \cite[I.3]{BoutotCarayol}. This links the objects studied here with the Bruhat-Tits building and the Berkovich embedding of the Bruhat-Tits building into the analytification of the flag variety. The present paper was motivated by this connection and the study done in \cite{PSS}.
\end{rem}

\section{Logarithmic differential operators on $\bbX_n$}\label{diff_op_Xn}

We refer to \cite{Montagnon} for a systematic discussion of sheaves of logarithmic differential operators. For $n \ge 1$ we equip $\bbX_n$ with the log structure defined by its normal crossings divisor $\{p=0\}$. However, here we will not use the theory as developed in \cite{Montagnon}, but rather work with a more elementary approach.

\subsection{The logarithmic tangent sheaf on $\bbX_n$}\label{log_tangent_sheaf}

\begin{para} For the purposes of this paper we consider the sheaf $\cD_{\bbX_n,\log}$ of logarithmic differential operators on $\bbX_n$ as being generated as a subsheaf of ${\mathcal End}_\Zp(\cO_{\bbX_n},\cO_{\bbX_n})$ by the logarithmic tangent sheaf $\cT_{\bbX_n,\log}$. (This is as in \cite[1.3]{Kato94}.) The restriction of $\cT_{\bbX_n,\log}$ to an open affine subset $\bbX^{(\nu)}_{\underline{a}}$, $\underline{a} \in \cR_\infty \times \cR^{\nu-1}$, cf. \ref{explicit_covering}, is generated by a differential operator $D$ (over $\Zp$) with the properties

$$D(x^{(\nu-1)}_{\underline{a}}) = x^{(\nu-1)}_{\underline{a}} \;, \hskip10pt D(z^{(\nu)}_{\underline{a}}) = -z^{(\nu)}_{\underline{a}} \;,$$

\vskip8pt

cf. \ref{local_eqn}. $D$ has the property that

$$D(x^{(\nu-1)}_{\underline{a}}z^{(\nu)}_{\underline{a}}) = x^{(\nu-1)}_{\underline{a}}D(z^{(\nu)}_{\underline{a}})+z^{(\nu)}_{\underline{a}}D(x^{(\nu-1)}_{\underline{a}}) = 0 \;,$$

\vskip8pt

and hence $D(x^{(\nu-1)}_{\underline{a}}z^{(\nu)}_{\underline{a}}-p) = 0$, so that $D$ preserves the ideal generated  by $x^{(\nu-1)}_{\underline{a}}z^{(\nu)}_{\underline{a}}-p$. Intuitively we may write

$$D = x^{(\nu-1)}_{\underline{a}}\partial_{x^{(\nu-1)}_{\underline{a}}} = - z^{(\nu)}_{\underline{a}}\partial_{z^{(\nu)}_{\underline{a}}} \;.$$

\vskip8pt

To put it another way, we may say that $\cT_{\bbX_n,\log}$ is locally on an open subscheme $\bbX^{(\nu)}_{\underline{a}}$ generated by

$$x^{(\nu-1)}_{\underline{a}}\partial_{x^{(\nu-1)}_{\underline{a}}}\, \mbox{ and } z^{(\nu)}_{\underline{a}}\partial_{z^{(\nu)}_{\underline{a}}} \, \;, \mbox{ with the relation } \;\;  x^{(\nu-1)}_{\underline{a}}\partial_{x^{(\nu-1)}_{\underline{a}}} = - z^{(\nu)}_{\underline{a}}\partial_{z^{(\nu)}_{\underline{a}}}  \,.$$

\vskip8pt

Denote by

$$pr_n: \bbX_n \lra \bbX_0 = \bbX$$

\vskip8pt

the canonical projection. Write, as in sec. \ref{smooth}, $\bbX_0 = \bbX = U_x \cup U_y$, where $U_x = \Spec(\Zp[x])$, $U_y = \Spec(\Zp[y])$, with $x$ and $y$ satisfying the relation $xy=1$. Let $\cI_{n,d} \sub \cO_\bbX$ be the ideal sheaf which is on $U_x$ associated to the ideal

$$\bigcap_{a \in \Zp/(p^n)} \left(x-a, p^n \right)^d \sub \Zp[x] = \cO_\bbX(U_x) \;,$$

\vskip8pt

and on $U_y$ associated to the ideal

$$\bigcap_{a \in \Zp/(p^n)} \left(y-a, p^n \right)^d \sub \Zp[y] = \cO_\bbX(U_y) \;.$$

\vskip8pt
\end{para}

Obviously, $\cI_{0,d} = \cO_\bbX$ for all $d$. In the following proposition, if $n=0$, we put $\cT_{\bbX_0,\log} = \cT_\bbX$.

\vskip8pt

\begin{prop}\label{direct_image} (a) $\cT_{\bbX_n,\log}$ is a subsheaf of the invertible sheaf $pr_n^*(\cT_\bbX)$.

\vskip8pt

(b) $(pr_n)_*(\cO_{\bbX_n}) = \cO_\bbX$.

\vskip8pt

(c) For all $n, d \ge 0$ one has $(pr_n)_*(\cT_{\bbX_n,\log}^{\otimes d}) = \cI_{n,d} \cT^{\otimes d}_\bbX$.
\end{prop}

\Pf (a) In order to see this we express the coordinate $x$ by the local coordinates $x^{(\nu-1)}_{\underline{a}}$ introduced in sec. \ref{covering}, and deduce a corresponding relation for $\partial_x$ and $\partial_{x^{(\nu-1)}_{\underline{a}}}$. (By symmetry it suffices to consider $x$.) To be precise, fix $1 \le \nu \le n$, $\underline{a} = (a_0, \ldots, a_{\nu-1}) \in \cR_\infty \times \cR^{\nu-1}$, and consider the open subset $\bbX^{(\nu)}_{\underline{a}} \sub \bbX_n$, cf. \ref{bbX_n}. Without loss of generality we may assume $a_0 \neq \infty$. Then we have $x-a = p^{\nu-1}x^{(\nu-1)}_{\underline{a}}$ where $a = a_0 + \ldots + a_{\nu-1}p^{\nu-1}$, cf. \ref{local_x_coord}. Hence

\begin{numequation}\label{relation_partials}
\partial_{x^{(\nu-1)}_{\underline{a}}} = p^{\nu-1}\partial_x \;, \; \mbox{and thus} \hskip8pt  x^{(\nu-1)}_{\underline{a}}\partial_{x^{(\nu-1)}_{\underline{a}}} = p^{\nu-1}x^{(\nu-1)}_{\underline{a}}\partial_x = (x-a)\partial_x \;.
\end{numequation}

This proves the assertion.

\vskip8pt

(b) The morphism $pr_n: \bbX_n \ra \bbX_0$ is a birational projective morphism of noetherian integral schemes, and $\bbX_0$ is normal. The assertion then follows exactly as in the proof of Zariski's Main Theorem as given in \cite[ch. III, Cor. 11.4]{HartshorneA}.

\vskip8pt

(c) {\it 1. The inclusion $(pr_n)_*(\cT_{\bbX_n,\log}^{\otimes d}) \sub \cI_{n,d} \cT^{\otimes d}_\bbX$.} Put $\bbX_n' = \bbX_n - pr_n^{-1}(\bbX(\Fp))$. This scheme is smooth over $\Zp$. The restriction of $pr_n$ induces an isomorphism

$$\bbX_n' \; \stackrel{\simeq}{\lra} \;  \bbX' = \bbX - \bbX(\Fp) \;,$$

\vskip8pt

and the restriction of $\cT_{\bbX_n,\log}$ to $\bbX_n'$ is the relative tangent sheaf of $\bbX_n'$ over $\Zp$ whose direct image under $pr_n$ is the relative tangent sheaf of $\bbX'$ over $\Zp$. Therefore, in order to understand $(pr_n)_*(\cT_{\bbX_n,\log}^{\otimes d})$ we need to investigate the stalks of this sheaf at the points in $\bbX(\Fp)$. We consider the point $P_0$ in $U_x = \Spec(\Zp[x]) \sub \bbX$ corresponding to the ideal $(x-a_0,p)$. Our aim is to understand the stalk of $(pr_n)_*(\cT_{\bbX_n,\log}^{\otimes d})$ at $P_0$.

\vskip8pt

By (a) we can consider the stalk of $(pr_n)_*(\cT_{\bbX_n,\log}^{\otimes d})$ at $P_0$ as a $\cO_{\bbX,P_0}$-submodule of the stalk of $\cT_\bbX^{\otimes d}$ at $P_0$. We consider thus an element

\begin{numequation}\label{function}
D = f(x)\partial_x^{\otimes d} \in \left(\cT_\bbX^{\otimes d}\right)_{P_0} \;,
\end{numequation}

$f(x) \in \cO_{\bbX,P_0} = \Zp[x-a_0]_{(x-a_0,p)}$, and want to find necessary and sufficient conditions for this element to be in the stalk of $(pr_n)_*(\cT_{\bbX_n,\log}^{\otimes d})$ at $P_0$. To this end, consider an open subset of $\bbX_n$ of the form

$$\bbX^{(1)}_{a_0} \cup \bbX^{(2)}_{a_0,a_1} \cup \ldots \cup \bbX^{(n)}_{a_0, \ldots, a_{n-1}} \cup \bbD^{(n)}_{a_0, \ldots, a_n} \;,$$

\vskip8pt

for a sequence $\underline{a} = (a_0, \ldots, a_n) \in \cR^{n+1}$. Consider the local coordinate $x^{(n)}_{a_0, \ldots, a_{n-1}, a_n}$ on $\bbD^{(n)}_{a_0, \ldots, a_n}$ which we denote henceforth by $x^{(n)}$. Put $a = a_0+a_1p+\ldots+a_{n-1}p^{n-1}+ a_np^n$. The equation \ref{local_x_coord} shows that

\begin{numequation}\label{relation_partials_Dn}
x^{(n)} = \frac{1}{p^n}\left(x-a\right) \;, \; \mbox{hence} \hskip10pt \partial_{x^{(n)}} = p^n\partial_x \;, \; \mbox{and thus} \hskip10pt x^{(n)}\partial_{x^{(n)}} = \left(x-a\right)\partial_{x-a} \;. 
\end{numequation}

If $D$ is in the stalk of
$(pr_n)_*(\cT_{\bbX_n,\log}^{\otimes d})$ at $P_0$ then $D$ extends to the stalk of $\cT_{\bbX_n,\log}^{\otimes d}$ at the point $P_n \in \bbD^{(n)}_{a_0, \ldots, a_n}$ corresponding to the ideal $(x^{(n)},p)$. Therefore, $D$ can be written as

\begin{numequation}\label{key}
g(x^{(n)})\partial_{x^{(n)}}^{\otimes d} \;,
\end{numequation}

with a function $g(x^{(n)}) \in \cO_{\bbX_n,P_n} = \Zp[x^{(n)}]_{(x^{(n)},p)}$. Completing this latter ring with respect to its maximal ideal gives $\Zp[[x^{(n)}]]$, and so we can consider $g(x^{(n)}) = \sum_{k\ge 0} c_k (x^{(n)})^k$ as an element in $\Zp[[x^{(n)}]]$. Now we write \ref{key} as

$$g\left(\frac{1}{p^n}(x-a)\right) p^{nd}\partial_x^{\otimes d} \;.$$

\vskip8pt

Using the power series expansion for $g(x^{(n)})$ gives

$$g\left(\frac{1}{p^n}(x-a)\right) p^{nd} = \sum_{k \ge 0} c_kp^{-nk+nd}(x-a)^k \;.$$

\vskip8pt

For $k \le d$ we have $p^{n(d-k)}(x-a)^k \in (x-a,p^n)^d$. And for $k > d$ we must have $c_kp^{-nk+nd} \in \Zp$ and so $c_kp^{-nk+nd}(x-a)^k$ is in $(x-a,p^n)^d$ too. The function $f(x)$ in \ref{function} is then contained in the ideal $(x-a,p^n)^d$ for all $a = a_0 + \ldots + a_{n-1}p^{n-1}$.
Hence we see that the stalk of
$(pr_n)_*(\cT_{\bbX_n,\log}^{\otimes d})$ at $P_0$ is contained in the stalk of $\cI_{n,d}\cT_\bbX^{\otimes d}$ at $P_0$. This is then true for all $\Fp$-rational points of $\bbX$. For the point at infinity one uses the equation \ref{local_y_coord}.

\vskip8pt

{\it 2. The inclusion $(pr_n)_*(\cT_{\bbX_n,\log}^{\otimes d}) \supset \cI_{n,d} \cT^{\otimes d}_\bbX$.} As above, we consider the point $P_0$ corresponding to the ideal $(x-a_0,p) \sub \Zp[x] = \cO_\bbX(U_x)$. For $1 \le \nu \le n$ consider an open affine subset $\bbX^{(\nu)}_{\underline{a}}$ of $\bbX_n$, as introduced in \ref{bbX_n} (cf. also \ref{explicit_covering}), where $\underline{a} = (a_0,a_1, \ldots, a_{\nu-1}) \in \cR_\infty \times \cR^{\nu-1}$. On $\bbX^{(\nu)}_{\underline{a}}$ we have the coordinate function $x^{(\nu-1)}_{\underline{a}}$, cf. \ref{local_eqn}, which is related to $x$ by

$$x = a_0+ \ldots +a_{\nu-1}p^{\nu-1} + p^{\nu-1}x^{(\nu-1)}_{\underline{a}} \;, \;\; \mbox{i.e.}\;, \hskip8pt x-a = p^{\nu-1}x^{(\nu-1)}_{\underline{a}} \;,$$

\vskip8pt

cf. \ref{local_x_coord}, where $a = a_0 + \ldots + a_{\nu-1}p^{\nu-1}$. Suppose $0 \le k \le d$ and consider the differential operator

$$D = p^{n(d-k)}(x-a)^k\partial_x^{\otimes d} \in (x-a,p^n)^d(\cT_\bbX^{\otimes d})_{P_0} \;.$$

\vskip8pt

We have $\partial_x = \frac{1}{p^{\nu-1}}\partial_{x^{(\nu-1)}_{\underline{a}}}$ and thus

\begin{numequation}\label{rewrite_D}
\begin{array}{rcl}
D & = & p^{n(d-k)}p^{k(\nu-1)-d(\nu-1)}(x^{(\nu-1)}_{\underline{a}})^k(\partial_{x^{(\nu-1)}_{\underline{a}}})^{\otimes d}  \\
&&\\
& = & p^{(n-\nu+1)(d-k)}(x^{(\nu-1)}_{\underline{a}})^k(\partial_{x^{(\nu-1)}_{\underline{a}}})^{\otimes d} \\
&&\\
&=& (z^{(\nu)}_{\underline{a}})^{(n-\nu+1)(d-k)}(x^{(\nu-1)}_{\underline{a}})^{(n-\nu+1)(d-k)}(x^{(\nu-1)}_{\underline{a}})^k(\partial_{x^{(\nu-1)}_{\underline{a}}})^{\otimes d}\\
&&\\
&=& (z^{(\nu)}_{\underline{a}})^{(n-\nu+1)(d-k)}(x^{(\nu-1)}_{\underline{a}})^{(n-\nu)(d-k)}(x^{(\nu-1)}_{\underline{a}})^d(\partial_{x^{(\nu-1)}_{\underline{a}}})^{\otimes d} \;.
\end{array}
\end{numequation}

Because of the term $(x^{(\nu-1)}_{\underline{a}})^d(\partial_{x^{(\nu-1)}_{\underline{a}}})^{\otimes d}$ on the last line of \ref{rewrite_D}, this shows that $D$ extends to  $\bbX^{(\nu)}_{\underline{a}}$. Here we have used the equation $z^{(\nu)}_{\underline{a}}x^{(\nu-1)}_{\underline{a}} = p$, cf. \ref{local_eqn}. To see that $D$ also extends to $\bbD^{(n)}_{\underline{b}}$, where here $\underline{b} = (a_0, a_1, \ldots, a_n)$, we use the coordinate $x^{(n)}_{\underline{b}}$ on $\bbD^{(n)}_{\underline{b}}$. The equations \ref{relation_partials_Dn} give then

$$D = p^{n(d-k)}p^{kn-dn}(x^{(n)}_{\underline{b}})^k(\partial_{x^{(n)}_{\underline{b}}})^{\otimes d}  = (x^{(n)}_{\underline{b}})^k(\partial_{x^{(n)}_{\underline{b}}})^{\otimes d} \;,$$

\vskip8pt

and this shows that $D$ extends to $\bbD^{(n)}_{\underline{b}}$.
If, more generally, we consider an element of the form $f(x)D$, where $f(x) \in \Zp[x]_{(x-a_0,p)}$ and $D$ is as before, then this will extend to a neighborhood of the special fiber of $\bbX^{(\nu)}_{\underline{a}}$ and $\bbD^{(n)}_{\underline{b}}$, respectively. \qed

\vskip8pt

\begin{cor}\label{direct_image_m} For all $n, d, m \ge 0$ one has

$$(pr_n)_*\left((\cT_{\bbX_n,\log}^{\otimes d})^{(m)}\right) = \cI_{n,d} (\cT^{\otimes d}_\bbX)^{(m)} = \frac{q^{(m)}_d!}{d!} \cI_{n,d} \cT^{\otimes d}_\bbX \;.$$

\vskip8pt

\end{cor}

\Pf The sheaf $\cT_{\bbX_n,\log}^{\otimes d}$ is a line bundle and the same reasoning as in the proof of \ref{trivialization} applies,
i.e., $(\cT_{\bbX_n,\log}^{\otimes d})^{(m)} = \frac{q^{(m)}_d!}{d!} \cT_{\bbX_n,\log}^{\otimes d}$. This equality is to be understood as in \ref{trivialization}. The statement then follows from \ref{direct_image}. \qed

\vskip8pt

Consider $\cI_{n,d} (\cT^{\otimes d}_\bbX)^{(m)}$ as a subsheaf of $(\cT^{\otimes d}_\bbX)^{(m)}$. The global sections of the former are thus contained in the global sections of the latter.

\vskip8pt

\begin{prop}\label{estimate} For all $n,d,m \ge 0$ one has the following inclusions

\begin{numequation}\label{inclusions} p^{nd}H^0\left(\bbX,(\cT^{\otimes d}_\bbX)^{(m)}\right) \;\; \sub \;\; H^0\left(\bbX,\cI_{n,d} (\cT^{\otimes d}_\bbX)^{(m)}\right) \;\; \sub \;\; p^{nc}H^0\left(\bbX,(\cT^{\otimes d}_\bbX)^{(m)}\right) \;,
\end{numequation}

as submodules of $H^0\left(\bbX,(\cT^{\otimes d}_\bbX)^{(m)}\right)$, where $c = \lceil d \frac{p-1}{p+1} \rceil$ is the smallest integer greater or equal to $d \frac{p-1}{p+1}$. In particular, for $d=1$ and any $n,m \ge 0$ we have

\begin{numequation}\label{estimate_tangent}
H^0\left(\bbX,\cI_{n,d} (\cT^{\otimes d}_\bbX)^{(m)}\right) \;\; = \;\; p^nH^0\left(\bbX,(\cT^{\otimes d}_\bbX)^{(m)}\right) \;.
\end{numequation}

\end{prop}

\Pf Because of \ref{trivialization} it suffices to treat the case $m=0$. By the very definition of $\cI_{n,d}$ one has $p^{nd}\cO_\bbX \sub \cI_{n,d}$ and thus $p^{nd} \cT^{\otimes d}_\bbX \sub \cI_{n,d} \cT^{\otimes d}_\bbX$. The inclusion on the left follows from this.
Furthermore, the statement is trivial for $n=0$ or $d=0$ (when $c=0$), and so we may assume that $n$ and $d$ are both positive.

\vskip8pt

To show the inclusion on the right we write global sections of $\cT^{\otimes d}_\bbX$ in the form $f(x)\partial_x^{\otimes d}$ with a polynomial $f(x) \in \Zp[x]$ of degree $\le 2d$.
Suppose $n \ge 1$ and $f(x)\partial_x^{\otimes d}$ is a global section
of $\cI_{n,d}\cT_\bbX^{\otimes d}$. Note that the reduction modulo $p$ of $\cI_{n,d}$ is an ideal sheaf on $\bbP^1_\Fp$ of degree $-(p+1)d$, which we denote by $\cI_{n,d,\Fp}$. Now, if $f$ is not divisible by $p$, then $(f \mod p) \partial_x^{\otimes d}$ would be a non-zero section of $\cI_{n,d,\Fp}\cT_{\bbP^1_{\Fp}}^{\otimes d}$ and this sheaf has degree $-(p+1)d+2d = (1-p)d <0$ (because $d>0$), hence a contradiction. Fix $a \in \Zp$ and write

\begin{numequation}\label{polynomial}
f(x) = \sum_{i=0}^d g_i(x)p^{ni}(x-a)^{d-i} \in (x-a,p^n)^d \;,
\end{numequation}

with polynomials $g_i(x) \in \Zp[x]$. We have seen that $f$ is divisible by $p$, hence $g_0$ is divisible by $p$.
Consider $\frac{1}{p}f(x) = \frac{g_0(x)}{p}(x-a)^d + \sum_{i=1}^d g_i(x)p^{ni-1}(x-a)^{d-i}$ and apply the previous reasoning. Doing this repeatedly shows that $g_0(x)$ is in fact divisible by $p^n$, and we find

$$f_1(x) \stackrel[\rm df]{}{=} \frac{1}{p^n}f(x) = \frac{g_0(x)}{p^n}(x-a)^d + \sum_{i=1}^d g_i(x)p^{n(i-1)}(x-a)^{(d-1)-(i-1)} \;,$$

\vskip8pt

and this polynomial is in $(x-a,p^n)^{d-1}$. This shows that $f_1(x) \partial_x^d$ is a global section of $\cI_{n,d-1}\cT_\bbX^{\otimes d}$. If $f_1$ is not divisible by $p$, then the same reasoning as above shows that $(f_1 \mod p) \partial_x^d$ gives rise to a non-zero global section of $\cI_{n,d-1,\Fp}\cT_{\bbP^1_{\Fp}}^{\otimes d}$ and this sheaf has degree $-(p+1)(d-1)+2d = (1-p)d + p+1$. If this number is negative we arrive at a contradiction. Suppose this number is non-negative. Arguing as above shows then that $f_1$ must be divisible by $p^n$, and hence $f$ is divisible by $p^{2n}$. Running the same arguments repeatedly proves that if $(1-p)d + j(p+1) <0$ we must have that $f$ is divisible by $p^{n(j+1)}$. Now the assertion follows because $c-1$ is the largest possible value for $j$. \qed

\vskip8pt

\begin{rem} The exponent $nc$ of $p$ on the right side of \ref{inclusions} is likely not the largest possible exponent for all $n$ and $d$. While it is interesting to find the largest possible exponent of $p$ for the inclusion on the right
side of \ref{inclusions}, the most optimistic guess that it be $nd$ is in general false. Consider for instance the case when $n=1$ and $d=p$. Then $p^{p-1}(x^p-x) \partial_x^{\otimes p}$ is a global section of $\cI_{1,p}\cT_\bbX^{\otimes p}$ as can be checked easily. We thus see that the optimal exponent would be at most $p-1$ and this is indeed equal to $\lceil p \frac{p-1}{p+1} \rceil$ for all $p$. Moreover, $p^{k(p-1)}(x^p-x)^k \partial_x^{\otimes kp}$ is a global section of $\cI_{1,kp}\cT_\bbX^{\otimes kp}$ for all $k$ and $p$, and we thus see that the exponent is at most $k(p-1) = kp \frac{p-1}{p}$. As a consequence, we see that the ratio $\frac{\mbox{{\scriptsize optimal exponent}}}{nd}$ is bounded by $\frac{p-1}{p}$ for $n=1$. Similar examples probably exist for arbitrary $n$.
\end{rem}

\vskip8pt

\subsection{Differential operators on $\bbX_n$ and distribution algebras}

Let $\cD_{\bbX_n}^{(m)} = \cD_{\bbX_n,\log}^{(m)}$ be the sheaf of logarithmic differential operators on $\bbX_n$ of level $m$. As an $\cO_{\bbX_n}$-module it is on an open affine subset $\bbX_{\underline{a}}^{(\nu)} \sub \bbX_n$, cf. \ref{explicit_covering}, locally generated by logarithmic differential operators

$$q^{(m)}_d! {D \choose d} \;$$

\vskip8pt

where

$$D = x^{(\nu-1)}_{\underline{a}}\partial_{x^{(\nu-1)}_{\underline{a}}} = - z^{(\nu)}_{\underline{a}}\partial_{z^{(\nu)}_{\underline{a}}} \;.$$

\vskip8pt

is a local section of the logarithmic tangent sheaf $\cT_{\bbX_n,\log}$, cf. \ref{log_tangent_sheaf}. On the open subscheme $\bbD^{(n)}_{\underline{b}}$ with coordinate function $x^{(n)}_{\underline{b}}$ it is generated by 

$$\frac{q^{(m)}_d!}{d!} \partial_{x^{(n)}_{\underline{b}}}^d \;.$$ 

\vskip8pt

Denote by $H^0(\bbX_n,\cD_{\bbX_n}^{(m)})^\wedge$ the $p$-adic completion of
$H^0(\bbX_n,\cD_{\bbX_n}^{(m)})$ and put $H^0(\bbX_n,\cD_{\bbX_n}^{(m)})^\wedge_\bbQ = H^0(\bbX_n,\cD_{\bbX_n}^{(m)})^\wedge \otimes_\bbZ \bbQ$.

\vskip8pt

\begin{thm}\label{thm-n} Given $n$ let $n' = \lfloor n\frac{p-1}{p+1} \rfloor$ be the greatest integer less or equal to $n\frac{p-1}{p+1}$. Then we have natural inclusions

$$\cD^{\rm an}(\bbG(n)^\circ)_{\theta_0} \; \hra \; \varinjlim_m H^0(\bbX_n,\cD_{\bbX_n}^{(m)})^\wedge_\bbQ \; \hra \; \cD^{\rm an}(\bbG(n')^\circ)_{\theta_0} \;.$$

\vskip8pt
\end{thm}

\Pf {\it 1. The inclusion on the left side.} The inclusion $\bbG(n)^\circ \sub \bbG(0)^\circ$ induces an embedding

$$\cD^{\rm an}(\bbG(n)^\circ)_{\theta_0} \hra \cD^{\rm an}(\bbG(0)^\circ)_{\theta_0} \;,$$

\vskip8pt

and the right hand side is canonically isomorphic to

$$\varinjlim_m H^0(\bbX,\cD_\bbX^{(m)})^\wedge_\bbQ \;,$$

\vskip8pt

by \ref{completion} and \ref{thm-0}. On the other hand, arguing as in the proof of \ref{direct_image}, part (a), one sees that $\cD_{\bbX_n}^{(m)}$ is naturally a subsheaf of $pr^*\left( \cD_\bbX^{(m)}\right)$, and so $H^0(\bbX_n,\cD_{\bbX_n}^{(m)}) \hra H^0(\bbX,\cD_\bbX^{(m)})$. The inclusion in question is thus understood to be an inclusion inside $\varinjlim_m H^0(\bbX,\cD_\bbX^{(m)})^\wedge_\bbQ$.

\vskip8pt

Now use \ref{dist_alg_G(n)} and the explicit form of the generators of $U(p^n\frg_\Zp)^{(m)}$ in \ref{generators_n}. Consider such an element

$$q^{(m)}_{\nu_1}! \frac{(p^ne)^{\nu_1}}{\nu_1!}  \cdot q^{(m)}_{\nu_2}! p^{n\nu_2}{h_1 \choose \nu_2} \cdot q^{(m)}_{\nu_3}! p^{n\nu_3} {h_2 \choose \nu_3} \cdot q^{(m)}_{\nu_4}! \frac{(p^nf)^{\nu_4}}{\nu_4!} \;.$$

\vskip8pt

Its image under the canonical map

$$\xi^{(m)}: U(\frg_\Zp)^{(m)} \lra H^0(\bbX,\cD^{(m)}_\bbX)$$

\vskip8pt

cf. \ref{map_U}, is

$$\frac{q^{(m)}_{\nu_1}!}{\nu_1!} (p^n \partial_x)^{\nu_1} \cdot \frac{q^{(m)}_{\nu_2}!}{\nu_2!} p^{n\nu_2}(-x)^{\nu_2} \partial_x^{\nu_2} \cdot \frac{q^{(m)}_{\nu_3}!}{\nu_3!} p^{n\nu_3} x^{\nu_3} \partial_x^{\nu_3} \cdot \frac{q^{(m)}_{\nu_4}!}{\nu_4!} (-p^nx^2\partial_x)^{\nu_4} \;.$$

\vskip8pt

The first and last term are of the form

$$\frac{q^{(m)}_\nu!}{\nu!} \left(p^n (\mbox{global section of } \cT_\bbX)\right)^\nu \;.$$

\vskip8pt

Because $H^0(\bbX_n,\cT_{\bbX_n,\log}) = p^nH^0(\bbX,\cT_\bbX)$, cf. \ref{estimate_tangent}, we see that these terms are in \linebreak $H^0(\bbX_n,\cD_{\bbX_n}^{(m)})$. For the second and third term we consider an open affine subset $\bbX^{(\mu)}_{\underline{a}}$. Let $x^{(\mu-1)} = x^{(\mu-1)}_{\underline{a}}$ be the coordinate on  $\bbX^{(\mu)}_{\underline{a}}$ as in \ref{local_eqn}.  Use \ref{relation_partials}, i.e., $\partial_{x^{(\mu-1)}} =  p^{\mu-1}\partial_x$, and write 

\begin{numequation}\label{calc}
\begin{array}{rcl}
\frac{q^{(m)}_\nu!}{\nu!} p^{n\nu}x^\nu\partial_x^\nu 
& = & \frac{q^{(m)}_\nu!}{\nu!} p^{n\nu}(x-a+a)^\nu\partial_x^\nu  \\
&&\\
& = & \sum_{k=0}^\nu \frac{q^{(m)}_\nu!}{(q^{(m)}_k!)(q^{(m)}_{\nu-k}!)} p^{nk} \cdot a^{\nu-k} \cdot  \frac{q^{(m)}_k!}{k!}(x-a)^k \partial_x^k \cdot \frac{q^{(m)}_{\nu-k}!}{(\nu-k)!}(p^n\partial_x)^{\nu-k} 
\end{array}
\end{numequation}

By what we have observed before we find that the term $\frac{q^{(m)}_{\nu-k}!}{(\nu-k)!}(p^n\partial_x)^{\nu-k}$ is a global section of $\cD_{\bbX_n}^{(m)}$. The relation $p^{\mu-1}x^{(\mu-1)} = x-a$ together with \ref{relation_partials} gives 

$$\frac{q^{(m)}_k!}{k!}(x-a)^k \partial_x^k = \frac{q^{(m)}_k!}{k!}(x^{(\mu-1)})^k \partial_{x^{(\mu-1)}}^k \;,$$

\vskip8pt

and so extends to a section of $\cD_{\bbX_n}^{(m)}$ over $\bbX^{(\mu)}_{\underline{a}}$. It is a straightforward exercise to see that $\frac{q^{(m)}_\nu!}{(q^{(m)}_k!)(q^{(m)}_{\nu-k}!)}$ is always an integer, and $\frac{q^{(m)}_\nu!}{\nu!} p^{n\nu}x^\nu\partial_x^\nu$ therefore extends to a section of $\cD_{\bbX_n}^{(m)}$ over $\bbX^{(\mu)}_{\underline{a}}$. Finally, we consider the subscheme $\bbD^{(n)}_{\underline{b}}$. Let $x^{(n)} = x^{(n)}_{\underline{b}}$ be the coordinate on $\bbD^{(n)}_{\underline{b}}$, as in \ref{bbX_n}, where $\underline{b} = (a_0, \ldots,a_n)$. Put $b = a_0+ \ldots+a_np^n$. Writing $x = (x-b)+b$, we can perform exactly the same calculation \ref{calc} as above, using \ref{relation_partials_Dn}, and find that $\frac{q^{(m)}_\nu!}{\nu!} p^{n\nu}x^\nu\partial_x^\nu$ extends to a section of $\cD_{\bbX_n}^{(m)}$ over $\bbD^{(n)}_{\underline{b}}$. 
So we can conclude that the terms $\frac{q^{(m)}_\nu!}{\nu!} p^{n\nu}x^\nu\partial_x^\nu$ are in $H^0(\bbX_n,\cD_{\bbX_n}^{(m)})$. 

\vskip8pt

The image of $\xi^{(m)}$ thus lies in $H^0(\bbX_n,\cD_{\bbX_n}^{(m)})$. Passing to the completions and the direct limit over $m$ we find that $\xi^{(m)}$ induces a map

$$D^{an}(\bbG(n)^\circ) \lra \varinjlim_m H^0(\bbX_n,\cD_{\bbX_n}^{(m)})^\wedge_\bbQ \;,$$

\vskip8pt

which makes the diagram

$$\xymatrix{
D^{an}(\bbG(n)^\circ) \ar[r] \ar[d] & \varinjlim_m H^0(\bbX_n,\cD_{\bbX_n}^{(m)})^\wedge_\bbQ \ar[d] \\
D^{an}(\bbG(0)^\circ)_{\theta_0} \ar[r] & \varinjlim_m H^0(\bbX,\cD_\bbX^{(m)})^\wedge_\bbQ
}$$

\vskip8pt

commute. The lower horizontal arrow is an isomorphism and the right vertical arrow is injective. The assertion now follows.

\vskip8pt

{\it 2. The inclusion on the right side.} For this inclusion consider the diagram

$$\xymatrix{
\varinjlim_m H^0(\bbX_n,\cD_{\bbX_n}^{(m)})^\wedge_\bbQ \ar@{-->}[r] \ar[d] & D^{an}(\bbG(n')^\circ)_{\theta_0}  \ar[d] \\
\varinjlim_m H^0(\bbX,\cD_\bbX^{(m)})^\wedge_\bbQ \ar[r] & D^{an}(\bbG(0)^\circ)_{\theta_0}
}$$

\vskip8pt

where the vertical arrows are injective and we have to show the existence of the dashed arrow. Let $N(m)$ be such that the cokernel of the canonical map

$$U(\frg_\Zp)^{(m)} \ra H^0(\bbX,\cD_\bbX^{(m)})$$

is annihilated by $p^{N(m)}$, cf. \ref{cokernel_bdd_torsion} (b).
Furthermore, consider the subsheaf $\cD_{\bbX_n,d}^{(m)}$ of
logarithmic differential operators of level $m$ and degree $\le
d$. Similarly, let $U(p^{n'}\frg_\Zp)^{(m)}_d$ be the submodule of
elements of degree $\le d$ as defined right before \ref{map}.
Then, in order to prove the existence of the dashed arrow in the
diagram above, it suffices to prove the existence of a map

$$p^{N(m)}H^0(\bbX_n,\cD_{\bbX_n,d}^{(m)}) \dashrightarrow U(p^{n'}\frg_\Zp)^{(m)}_d \;,$$

\vskip8pt

which makes the corresponding diagram

$$\xymatrix{
p^{N(m)}H^0(\bbX_n,\cD_{\bbX_n,d}^{(m)}) \ar@{-->}[r] \ar[d] &   U(p^{n'}\frg_\Zp)^{(m)}_d \ar[d] \\
p^{N(m)}H^0(\bbX,\cD_\bbX^{(m)}) \ar[r] & U(\frg_\Zp)^{(m)}_d
}$$

\vskip8pt

commute. We do this by induction over $d$. This is obvious for $d=0$. For the induction step we can pass to the corresponding graded object in degree $d$ and thus consider

$$\xymatrix{
p^{N(m)}H^0(\bbX_n,(\cT_{\bbX_n,\log}^{\otimes d})^{(m)}) \ar@{-->}[r] \ar[d] &   U(p^{n'}\frg_\Zp)^{(m)}_d/U(p^{n'}\frg_\Zp)^{(m)}_{d-1} \ar[d] \\
p^{N(m)}H^0(\bbX,(\cT_\bbX^{\otimes d})^{(m)}) \ar[r] & U(\frg_\Zp)^{(m)}_d/U(\frg_\Zp)^{(m)}_{d-1}
}$$

\vskip8pt

Note that

$$U(p^{n'}\frg_\Zp)^{(m)}_d/U(p^{n'}\frg_\Zp)^{(m)}_{d-1} \;\; = \;\; p^{dn'} \left(U(\frg_\Zp)^{(m)}_d/U(\frg_\Zp)^{(m)}_{d-1}\right) \;.$$

\vskip8pt

By \ref{direct_image} and \ref{estimate} we have an inclusion

$$H^0(\bbX_n,(\cT_{\bbX_n,\log}^{\otimes d})^{(m)}) \sub p^{nc(d)}H^0(\bbX,(\cT_\bbX^{\otimes d})^{(m)}) \;,$$

\vskip8pt

where $c(d) = \left\lceil d \frac{p-1}{p+1} \right\rceil$.
The assertion now follows from the following inequalities:

$$nc(d) = n \left\lceil d \frac{p-1}{p+1} \right\rceil \ge nd \frac{p-1}{p+1}  \ge d \left\lfloor n \frac{p-1}{p+1} \right\rfloor = dn' \;.$$

\qed

\vskip8pt

\begin{rem}\label{arithmetic} We recall that $\frX_n$ denotes the completion of $\bbX_n$ along its special fiber, and we let $\sD_{\frX_n}^{(m)} = \widehat{\cD}_{\bbX_n}^{(m)}$ be the $p$-adic completion of the sheaf $\cD_{\bbX_n}^{(m)}$. Consider these as sheaves on $\frX_n$. Put $\sD_{\frX_n,\bbQ}^{(m)} = \sD_{\frX_n}^{(m)} \otimes_\bbZ \bbQ$ and

$$\sD^\dagger_{\frX_n,\bbQ} \stackrel[\rm df]{}{=} \varinjlim_m \sD^{(m)}_{\frX_n,\bbQ} \;.$$

\vskip8pt

Then, as is not difficult to see, there is a canonical injective ring homomorphism

\begin{numequation}\label{compl_vs_arithmetic}
\varinjlim_m H^0(\bbX_n,\cD_{\bbX_n}^{(m)})^\wedge_\bbQ \; \hra \; H^0(\frX_n,\sD^\dagger_{\frX_n,\bbQ}) \;.
\end{numequation}

\vskip8pt

The same reasoning as in \cite[Prop. 3.2]{Huyghe97} shows that
this map is an isomorphism, if $H^1(\bbX_n,\cD_{\bbX_n}^{(m)})$ is
annihilated by some fixed power of $p$. This question in turn is
closely connected to the question whether $\frX_n$ is
$\sD^\dagger_{\frX_n,\bbQ}$-affine, a problem we plan to discuss
in a future paper.

\end{rem}

\bibliographystyle{plain}
\bibliography{mybib}

\end{document}